 \documentclass[preprint,review,12pt]{elsarticle}

\usepackage{amssymb}
 \usepackage{amsthm,amsmath}
\usepackage{mathrsfs}

 \numberwithin{equation}{section}

\newtheorem{thm}{Theorem}[section]
\newtheorem{lem}[thm]{Lemma}

\theoremstyle{definition}

\newcommand{\beq}{\begin{eqnarray}}
\newcommand{\eeq}{\end{eqnarray}}
\newcommand{\beqno}{\begin{eqnarray*}}
\newcommand{\eeqno}{\end{eqnarray*}}

\newenvironment {Proof of Theorem} {\noindent {\bf Proof of Theorem 2.1}}{\quad $\square$\par\vspace{3mm}}

\allowdisplaybreaks

\journal{Nonlinear of Analysis, Real World and Applications}

\begin{document}
\newcommand{\D}{\displaystyle}
\begin{frontmatter}



\title{Local existence of unique strong solution to non-isothermal model for incompressible nematic liquid crystals in 3D}


\author{Shijin Ding$^{1}$}
\author{Quanrong Li$^{1*}$\corref{cor1}}
\cortext[cor2]{Corresponding author. Email: 787177237@qq.com}
\address{$^{1}$School of Mathematical Sciences,
South China Normal University,  Guangzhou  510631, China.}

\begin{abstract}
In this paper, we consider the non-isothermal model for incompressible flow of nematic liquid crystals in three dimensions and prove the local existence and uniqueness of the strong solution with
periodic initial conditions on $ \mathbb{T}^3$.
\end{abstract}

\begin{keyword}
incompressible flow\sep liquid crystals\sep non-isothermal\sep local existence, uniqueness.\\
{\em AMS Subject Classification:} 76N10\sep 35Q30\sep 35R35.
\end{keyword}

\end{frontmatter}

\newpage
\section{Introduction}
Liquid crystals, which exist in an intermediate state between isotropic liquid and solid, are materials with rheological properties. They are first found in the  extractant of nerve fiber by Prussia doctors Rudolf Virchow, et al., in the middle of the nineteenth century, but not until the 1930s, when some scientists successively applied liquid crystals to temperature sensors, display elements and optical memory devices, did people pay wide
attention to these materials. With the development of the modern industry, liquid crystals have been widely applied to thin type display devices and become the well-known industrial materials.

Many different types of liquid crystals phases have been observed in practical experiment, distinguished by their special optical properties. When viewed under a microscope with a polarized light source, different types of liquid crystals will appear distinct textures, which means the liquid crystals molecules are oriented in different directions. However, the molecules are well ordered within a specific domain. Real-world applications motivates the theoretical research, for proper functioning of many practical devices depends on special optical properties of liquid crystals in the presence or absence of an electric field. As for describing the dynamic of liquid crystal flow by continuum theories, there are lots of attempts. Here we would like to refer to some celebrate papers, such as\cite{1,2}, where Ericksen and Leslie provided a description on mathematics to various properties of liquid crystals, while F. H. Lin and C. Liu \cite{3} for the first time to analysis the model mathematically. As to thermotropic liquid crystals, there are three main types approved by the general, termed \textit{sematic, nematic}, and \textit{cholesteric}.

In this paper, we consider a non-isothermal model of nematic liquid crystals, which was based on the spirit of the simplified version of the Ericksen-Leslie system and proposed by E. Feireisl et al. in \cite{4}. In this model, the time evolution of the velocity field $u$ is governed by the standard incompressible Navier-Stokes equations with a non-isotropic stress tensor depending on $\nabla u$ and $\nabla d$. In addition, the transport coefficients vary with temperature. The dynamic behavior of the direction field $d$ is described by a  penalized Ginzburg-Landau type system, while the heat equation of the absolute temperature $\theta$ is according to a variant of \textit{Fourier's low}, in which the direction field $d$ is taken into account. More precisely, we study the following systems:
\begin{equation}
\begin{cases}
\partial_t\rho+\mathrm{div}(\rho u)=0\\
\mathrm{div}(u)=0\\
\partial_t(\rho u)+\mathrm{div}(\rho u\otimes u)+\nabla p=\mathrm{div}(\mathbb{S}+\rho\sigma^{nd})+\rho f\\
\partial_t(\rho\theta)+\mathrm{div}(\rho \theta u)+\mathrm{div}q=(\mathbb{S}+\rho\sigma^{nd}):\nabla u\\
\partial_td+u\cdot\nabla d+\partial_dW(d)=\frac{1}{\rho}\mathrm{div}(\rho\nabla d),\\
\end{cases}
\end{equation}
where
\begin{align*}
&\mathbb{S}=\mu(\theta)(\nabla u+\nabla ^{T}u),\sigma^{nd}=-\lambda(\theta)\nabla d\odot\nabla d,\\
&q=-\kappa(\theta)\nabla\theta-(\kappa_{\parallel}-\kappa_{\perp})(\theta)(d\cdot\nabla\theta)d,\\
&\nabla d\odot\nabla d\ \text{is a}\ 3\times 3\ \text{matrix}, (\nabla d\odot\nabla d)_{ij}=\partial_{x_i}d\cdot\partial_{x_j}d,i,j=1,2,3.
\end{align*}
In addition,~$\rho $ represents the mass density,~$u$ represents the velocity field,~$\theta$ stands for the absolute temperature and ~$d$ to be the average direction field of the molecules of liquid crystals. Moreover, $\mu (\theta),\lambda(\theta) $ are viscosity coefficients and $f$ is the external force, while~$\kappa(\theta)$ and $(\kappa_{\parallel}-\kappa_{\perp})(\theta)$ are positive functions about the absolute temperature.

In \cite{4}, Feireisl et al. only investigated a special case of the above system, that is, let ~$\rho \equiv 1, f\equiv 0 $ and tried to solve the following initial boundary value problem:
\begin{equation}
\begin{cases}
\mathrm{div}u=0\\
\partial_tu+u\cdot\nabla u+\nabla p=\mathrm{div}(\mathbb{S}+\sigma^{nd})\\
\partial_t\theta+u\cdot\nabla\theta+\mathrm{div}q=(\mathbb{S}+\sigma^{nd}):\nabla u\\
\partial_td+u\cdot\nabla d+\partial_dW(d)=\Delta d,\\
\end{cases}
\end{equation}
\begin{align}
&(u,d,\theta)|_{t=0}=(u_{0},d_{0},\theta_{0}),\\
&u\cdot n|_{\partial\Omega}=0,[\mathbb{T} n]\times n|_{\partial\Omega}=0,\\
&q\cdot n|_{\partial\Omega}=0,\\
&\nabla_x d_i\cdot n|_{\partial\Omega}=0,\mathrm{for}\ i=1,2,3,
\end{align}
where $ W(d)=(|d|^{2} -1)^{2}, \mathbb{T}=\mathbb{S}+\sigma^{nd}-p\mathbb{I}$.

Assuming that the viscosity coefficients and initial data satisfy some natural conditions on regularity, Feireisl et al. obtained the global-in-time weak solution of the above initial boundary value problem in three dimensions by means of Galerkin approximation. In 2012, Feireisl et al. in \cite{5} introduced another non-isothermal model for nematic liquid crystal flow and similar result was proved.

Based on work of Feireisl et al., J. K. Li and Z. P. Xin \cite{6} also made some efforts on non-isothermal model for the dynamic of nematic liquid crystals flow. They reserved the constraint $|d|=1$, and modified some terms of the model that proposed by Feireisl et al. Precisely, they investigated the following system:
\begin{equation}
\begin{cases}
u_{t}+u\cdot\nabla u+\nabla p=\mathrm{div}(\mathbb{S}+\sigma^{nd}),\mathrm{div}u=0\\
d_{t}+u\cdot\nabla d=\Delta d+|\nabla d|^{2}d,|d|=1\\
\theta_{t}+u\cdot\nabla\theta=\Delta\theta+\mathbb{S}:\nabla u+|\Delta d+|\nabla d|^{2}d|^2,\\
\end{cases}
\end{equation}
where $\mathbb{S}=\mu(\theta)(\nabla u+\nabla^{T}u),\sigma^{nd}=-\nabla d\odot\nabla d.$

By means of $L^\infty$-type cut-off and Galerkin approximation, the global existence of weak solutions to the above initial boundary problem on $\mathbb{T}^2$ was obtained by assuming some suitable smoothness to the initial data.

The model we take into account comes from \cite{4}, which is proposed by Feireisl et al., and we will also do some simplification. That is, taking ~$\rho\equiv 1,f\equiv 0,\kappa(\theta)\equiv 1,(\kappa_\parallel-\kappa_\perp)(\theta)\equiv 0$ and finally we study the following problem:
\begin{align}
\begin{cases}
u_t+u\cdot\nabla u+\nabla p=\mathrm{div}(\mathbb{S}+\sigma^{nd})\\
\mathrm{div}u=0\\
\theta_t+u\cdot\nabla\theta=\Delta\theta+(\mathbb{S}+\sigma^{nd}):\nabla u\\
d_t+u\cdot\nabla d=\Delta d-(|d|^2-1)d,
\end{cases}
\end{align}
where $\mathbb{S}=\mu(\theta)(\nabla u+\nabla^Tu),\sigma^{nd}=-\lambda(\theta)\nabla d\odot\nabla d, $
subject to the initial data:
\begin{align}
&(u,\theta,d)|_{t=0}=(u_0,\theta_0,d_0),\nonumber\\
&\mathrm{div}u_0=0,\mathrm{ess}\inf_\Omega\theta_0=\underline{\theta}_0>0,|d_0|=1,a.e.\text{ÓÚ}\Omega=[-D,D]^3\subset\mathbb{R}^3
\end{align}
and boundary conditions£º
\begin{align}
&u(x+De_i,t)=u(x-De_i,t),\theta(x+De_i,t)=\theta(x-De_i,t),\nonumber\\
&d(x+De_i,t)=d(x-De_i,t),\forall x\in\partial\Omega,i=1,2,3.
\end{align}

Although such a model may seem rather naive from the point of view of real-world applications, it does capture the essential mathematical features of the problem. Moreover, it is compatible with all underlying physical principles, especially with \textit{the first and the second laws of thermodynamics}. What we want to achieve is to prove the existence of unique local-in-time strong solution without any essential restriction to the initial data.

Throughout this paper, we assume that $\mu(\cdot),\lambda(\cdot)\in C^2([0,\infty)) $ and satisfy the following conditions:
\begin{align}
&\lambda^\prime(0)>0,\lambda(0)=0,\lambda(\theta)\leq \bar{\lambda},
0\leq\lambda^\prime(\theta)\leq\overline{\lambda^\prime},|\lambda^{\prime\prime}(\theta)|\leq\overline{\lambda^{\prime\prime}}\quad\forall \theta\geq 0,\\
&\underline{\mu}\leq\mu(\theta)\leq\overline{\mu},|\mu^\prime(\theta)|\leq\overline{\mu^\prime},|\mu^{\prime\prime}(\theta)|\leq\overline{\mu^{\prime\prime}}\quad\forall \theta\geq 0.
\end{align}
Here $\bar{\lambda},\overline{\lambda^\prime},\overline{\lambda^{\prime\prime}},\underline{\mu},\overline{\mu},\overline{\mu^\prime},\overline{\mu^{\prime\prime}}>0$ are constants depending on the properties of materials.

Before continuing, we would like to have some more words here on the construction of this system. When the direction field $d\equiv 1$, (1.8) is just the non-isentropic incompressible Navier-Stokes system with temperature-dependent viscosity coefficients and invariable density, the existence of global-in-time weak solutions to which had been proved by Feireisl et al. in \cite{7}. Furthermore, Y. Cho and H. Kim in \cite{8} discussed a more general non-isentropic incompressible Navier-Stokes system, where the transport coefficients depend on both temperature and density, and the existence of local-in-time unique strong solution is obtained, see also \cite{9}. However, for compressible cases, there are only one-dimensional relevant results at present, for which readers can refer to \cite{10,11,12,13}.

 On the other hand, when the absolute temperature is not taken into account, the viscosity coefficients will turn into constants and system (1.8) will become that one proposed by F. H. Lin and C. Liu in \cite{3}. In 1995, Lin and Liu in \cite{3} proved the existence of global-in-time weak solutions and unique classical solution under the assumption that the viscosity is a large enough constant, while in 2000 in \cite{14}, the existence of global weak solutions, local classical unique solution and global classical unique solution with large viscosity are obtained for a more general Ericksen-Leslie model with penalization. For non-penalized model, S. J. Ding and H. Y. Wen \cite{19} obtained the local existence and uniqueness of the strong solution in dimension two or three with vacuum, while without vacuum, they obtained the global existence and uniqueness of the solution with small initial data in dimension two and S. J. Ding et al.\cite{20} dealt with the three dimensional case. However, the existence of global weak solutions in three dimensions is now still an open problem in general, although the two dimensional problem was considered by F. H. Lin et al.\cite{21}. The most general idea is taking $\varepsilon\rightarrow 0$ and considering the limitation of the penalized model, which was first tried by Lin and Liu in \cite{15} for the isentropic case. For the recent motivation, we refer the readers to \cite{16,17}.

The rest of this paper is arranged as follows: in section 2, some notations, which will be used throughout this paper, and the main theorem are stated; in section 3, we carry out the Galerkin approximation and obtain the approximate solutions for fixed $m$, while in section 4, uniform {\it a priori} estimation of the approximate solutions for at least a finite time will be established, where, in order to overcome the strongly coupled properties and nonlinearity, we apply the Gronwall inequality to obtain the high order energy estimates for the approximate solutions; finally in section 5 and section 6, the existence and uniqueness of the local-in-time solution of system (1.8) are proved respectively.
\section{Notations and main theorem}
We introduce some notations used throughout this paper.
\begin{align*}
&\bullet V=\{v\in C^\infty_0(\Omega)|\mathrm{div}u=0\}, W=\text{closure of }V\ \text{in}\ H^1_0(\Omega).\\
&\bullet\int f:=\int_{\Omega}f\mathrm{d}x,L^q:=L^q,H^k:=H^k(\Omega)=W^{k,2}(\Omega),H^1_0:=H^1_0(\Omega).\\
&\bullet A\lesssim B\ \text{denotes that}\ A\leq CB,\ \text{where}\ C > 0\ \text{is an absolute positive constant}.\\
&\bullet\text{For}\ 3\times 3\ \text{matrixes}\ A=(a_{ij})_{1\leq i,j\leq 3}\ \text{and}\ B=(b_{ij})_{1\leq i,j\leq 3},\text{we denote}\\
&\quad A:B=\sum^3_{i,j=1}a_{ij}b_{ij}.
\end{align*}
In this paper, \textit{the Einstein summation convention} is used, that is, repeated index contains the sum
within an monomial.

The main result of the present paper is the following theorem:
\begin{thm}
Assume that ~$(u_0,\theta_0,d_0)\in H^2(\Omega)\cap W\times H^2\cap L^1\times H^3\cap L^\infty$. Then there exists $T^*>0$, such that for any positive time $T \in[0,T^*)$, there is a unique solution $(u,\theta,d,p)$ in $[0,T]$ solve the system (1.8), satisfying
\begin{align}
&u\in L^\infty(0,T;H^2\cap W)\cap L^2(0,T;H^3),u_t\in L^2(0,T;H^1);\\
&d\in L^\infty(0,T;H^3)\cap L^2(0,T;H^4)\cap L^\infty([0,T]\times\Omega),d_t\in L^2(0,T;H^2);\\
&\theta\in L^\infty(0,T;H^2)\cap L^2(0,T;H^3),\theta_t\in L^2(0,T;H^1);\\
&p\in L^\infty(0,T;H^1)\cap L^2(0,T;H^2).
\end{align}
\end{thm}
\section{Galerkin approximation}
In this section, we will get an approximate solution for the proposed problem by means of Galerkin approximation.

Denote
\begin{equation*}
X=\{u\in H^1_0\cap H^2|\nabla\cdot u=0\}
\end{equation*}
and its finite-dimensional subspaces
\begin{equation*}
X^m=\mathrm{span}\{\Phi^1,\Phi^2,\cdots,\Phi^m\},m=1,2,\cdots
\end{equation*}
where\ $\Phi^m,m=1,2,\cdots $~is an complete normal orthogonal base of~$X$, which means that ~$(\Phi^i,\Phi^j)=\delta_{ij}$,where~$\delta_{ij}=1$~if~$i=j$,~and~$\delta_{ij}=0$~if~$i\neq j$.

For fixed~$m$, and~$k=1,2,\cdots,m$, we study the following approximate system:
\begin{align}
\begin{cases}
(u^m_t,\Phi^k)+(u^m\cdot\nabla u^m,\Phi^k)+(\mathbb{S}^m,\nabla \Phi^k)=-(\sigma^{nd,m},\nabla\Phi^k)\\
\mathrm{div}u^m=0\\
\theta^m_t+u^m\cdot\nabla\theta^m=\Delta\theta^m+(\mathbb{S}^m+\sigma^{nd,m}):\nabla u^m\\
d^m_t+u^m\cdot\nabla d^m=\Delta d^m-(|d^m|^2-1)d^m\\
(u^m,\theta^m,d^m)|_{t=0}=(u^m_0,\theta_0, d_0)\\
u^m(x+De_i,t)=u^m(x-De_i,t),d^m(x+De_i,t)=d^m(x-De_i,t),\\
\theta^m(x+De_i,t)=\theta^m(x-De_i,t),\forall x\in\partial\Omega,i=1,2,3,\\
\end{cases}
\end{align}
where~$\mathbb{S}^m=\mu(\theta^m)(\nabla u^m+\nabla^Tu^m)$, $\sigma^{nd,m}=-\lambda(\theta^m)\nabla d^m\odot\nabla d^m $, $u_0^m=\sum^m_{i=1}(u_0,\Phi^i)\Phi^i$,\\
$(f,g):=\int_{\Omega}f\cdot g\mathrm{d}x$\ is the inner product on $L^2(\Omega;\mathbb{R}^3)$.

We will try to find a solution to the above initial boundary value problem, in which the velocity field is  assumed to be $u^m(x,t)=\sum^m_{i=1}\tilde{g}^m_i(t)\Phi^i(x)$. Here~$\tilde{g}^m_i(t)$ satisfies~$\tilde{g}^m_i(0)=(u_0,\Phi^i),\tilde{g}^m_i(t)\in C^1[0,T],i=1,2,\cdots,m$.

\textit{Step 1:} Define a finite dimensional function space as
\begin{equation}
\mathcal{C}=\{g^m(t)=(g^m_1(t),g^m_2(t),\cdots,g^m_m(t))|g^m_i(t)\in C[0,T]\}
\end{equation}
and a closed convex subset of this space as
\begin{equation}
\mathcal{C}_K=\{g^m(t)\in(C[0,T])^m||g^m(t)|^2=\sum^m_{i=1}|g^m_i(t)|^2\leq K\}
\end{equation}
where $K, T$~are positive constants to be determined.

Then, in order to difine a map from $\mathcal{C}_K$ to itself, we proceed as follows. For any given function
\begin{equation}
v^m(x,t)=\sum_{i=1}^mg^m_i(t)\Phi^i(x),
\end{equation}
such that
\begin{equation}
g^m(t)=(g^m_1(t),g^m_2(t),\cdots,g^m_m(t))\in\mathcal{C}_K,
\end{equation}
we solve the equation
\begin{equation}
d_t^m+v^m\cdot\nabla d^m=\Delta d^m-(|d^m|^2-1)d^m
\end{equation}
 with the same initial condition as before, where~$v^m(x,t)$~is given in \textit{Step 1}.

\textit{Step 2:}
We should note that, this is a semi-linear parabolic partial differential equation. By the classical theory of parabolic partial differential equations and fixed point theorem, for any $T>0$, there is a unique function~$d^m\in L^\infty(0,T;H^3(\Omega))$ solves the following initial boundary problem:
\begin{equation}
\begin{cases}
d^m_t+v^m\cdot\nabla d^m=\Delta d^m-(|d^m|^2-1)d^m\\
d^m|_{t=0}=d_0\in H^3(\Omega, S^2)\\
d^m(x+De_i,t)=d^m(x-De_i,t),\forall x\in\partial\Omega,i=1,2,3.
\end{cases}
\end{equation}
Therefore, we have defined a continuous map~$\mathcal{T}_1:v^m\mapsto d^m$.

\textit{Step 3}: Solve equation
\begin{align}
\theta^m_t+v^m\cdot\nabla\theta^m=\Delta\theta^m+[\mu(\theta^m)(\nabla v^m+\nabla^Tv^m)-\lambda(\theta^m)\nabla d^m\odot\nabla d^m]:\nabla v^m,
\end{align}
where~$v^m(x,t)$~is the approximate function from \textit{Step 1} and~$d^m(x,t)$~from \textit{Step 2}.

Similar to \textit{Step 2}, for any $T>0$, there is a unique function~$\theta^m\in L^\infty(0,T;H^2(\Omega))$ solves the following initial boundary problem:
\begin{align}
\begin{cases}
\theta^m_t+v^m\cdot\nabla\theta^m=\Delta\theta^m+\frac{\mu(\theta^m)}{2}|\nabla v^m+\nabla^Tv^m|^2-\lambda(\theta^m)\nabla d^m\odot\nabla d^m:\nabla v^m\\
\theta^m|_{t=0}=\theta_0\in H^2(\Omega)\\
\theta^m(x+De_i,t)=\theta^m(x-De_i,t),\forall x\in\partial\Omega,i=1,2,3.
\end{cases}
\end{align}
Here we define a continuous map~$\mathcal{T}_2:(v^m,d^m)\mapsto\theta^m$.

In addition, according to the differential mean value theorem, $\lambda(\theta)=\lambda(\theta)-\lambda(0)=\lambda^\prime(\xi\theta)\theta,0<\xi<1$. On the other hand, since the initial data satisfies that~$\mathrm{ess}\inf_\Omega\theta_0=\underline{\theta}_0$, by maximum principle, we have $\mathrm{ess}\inf_\Omega\theta(x,t)\geq\underline{\theta}_0>0,a.e\ \Omega\times[0,T]$.

\textit{Step 4}: Solve equation
\begin{align}
(u^m_t,\Phi^k)+(u^m\cdot\nabla u^m,\Phi^k)=-(\mathbb{S}^m+\sigma^{nd,m},\nabla\Phi^k),k=1,2,\cdots,m,
\end{align}
where~$\mathbb{S}^m=\mu(\theta^m)(\nabla u^m+\nabla^Tu^m), \sigma^{nd,m}=-\lambda(\theta^m)\nabla d^m\odot\nabla d^m$ and $d^m(x,t),\theta^m(x,t)$ are solutions from $\textit{Step 2}$ and $\textit{Step 3}$ respectively.

Substituting~$u^m(x,t)=\sum^m_{i=1}\tilde{g}^m_i(t)\Phi^i(x)$~into the equation above, by the properties of ~$\Phi^i(i=1,\cdots,m)$ we obtain an equivalent system of ordinary equations as follows:
\begin{align}
\frac{\mathrm{d}}{\mathrm{d}t}\tilde{g}^m_k(t)=-\sum^m_{i=1}A_i^k(t)\tilde{g}^m_i(t)-\sum^m_{i,j=1}e^k_{ij}\tilde{g}^m_i(t)\tilde{g}^m_j(t)+ f^k(t),
\end{align}
where
\begin{align*}
&A_i^k(t)=\int\mu(\theta^m)(\nabla \Phi^i+\nabla^T\Phi^i):\nabla\Phi^k,\\
&e^k_{ij}=\int\Phi^i\cdot\nabla\Phi^j\cdot\Phi^k,\\
&f^k(t)=\int\lambda(\theta^m)\nabla d^m\odot\nabla d^m:\nabla\Phi^k,k=1,2,\cdots,m.
\end{align*}

Then, to determine the solution $u^m$, one finds that it suffices to solve the Cauchy problem:
\begin{align}
\begin{cases}
\frac{\mathrm{d}}{\mathrm{d}t}\tilde{g}(t)=G(\tilde{g},t)\\
\tilde{g}^m_k(0)=(u_0,\Phi^k),k=1,2,\cdots,m.
\end{cases}
\end{align}
Again we should note that
\begin{align}
&G^k(\tilde{g},t)=-\sum^m_{i=1}A_i^k(t)\tilde{g}^m_i(t)-\sum^m_{i,j=1}e^k_{ij}\tilde{g}^m_i(t)\tilde{g}^m_j(t)+ f^k(t),\nonumber\\
&\tilde{g}(t)=(\tilde{g}^m_1(t),\tilde{g}^m_2(t),\cdots,\tilde{g}^m_m(t)).\nonumber
\end{align}

By the theory of ordinary differential equations, this Cauchy problem admits  unique solution $\tilde{g}(t)$ in some finite time interval $[0,T_0]\subset[0,T]$~such that $\tilde{g}(t)\in(C^1[0,T_0])^m$, which yields $u^m=\sum^m_{i=1}\tilde{g}^m_i(t)\Phi^i(x)$ solves system (3.10).

This allows us to define the third continuous map~$\mathcal{T}_3:(d^m,\theta^m)\mapsto u^m$.

\textit{Step 5}: To sum up, we have defined a continuous map~$\mathcal{T}:v^m\mapsto u^m$~by the composition of the  maps $\mathcal{T}_1,\mathcal{T}_2$~and~$\mathcal{T}_3$. Furthermore, this map is equivalent to the following one:
\begin{align}
\widetilde{\mathcal{T}}:g(t)=(g^m_1(t),g^m_2(t),\cdots,g^m_m(t))\mapsto \tilde{g}(t)=(\tilde{g}^m_1(t),\tilde{g}^m_2(t),\cdots,\tilde{g}^m_m(t)),\nonumber
\end{align}
where~$v^m=\sum^m_{i=1}g^m_i(t)\Phi^i(x),u^m=\sum^m_{i=1}\tilde{g}^m_i(t)\Phi^i(x)$.

Next, we want to verify that the map~$\widetilde{\mathcal{T}}$~does satisfy conditions of Schauder fixed point theorem so that we can infer that system (3.1) has a solution~$(u^m,\theta^m,d^m)$.

The continuity of $\widetilde{\mathcal{T}}$ can be easily obtained by that of maps defined from \textit{Step 2} to \textit{Step 4}. In order to use Schauder fixed point theorem, we still need to prove that~$\tilde{g}(t)\in\mathcal{C}_K$~and~$\widetilde{\mathcal{T}}$ is a compact operator.

Actually, multiplying (3.11) by $\tilde{g}^m_k(t)$ and summing up $k$ from~$1$~to~$m$, we have
\begin{align}
\frac{1}{2}\frac{\mathrm{d}}{\mathrm{d}t}\sum^m_{k=1}(\tilde{g}^m_k(t))^2=-\sum^m_{i,k=1}A_i^k(t)\tilde{g}^m_i(t)\tilde{g}^m_k(t)-\sum^m_{k,i,j=1}e^k_{ij}\tilde{g}^m_i(t)\tilde{g}^m_j(t)\tilde{g}^m_k(t)
+\sum^m_{k=1} f^k(t)\tilde{g}^m_k(t).
\end{align}
Integrating by parts, we note that
\begin{align*}
e^k_{ij}=(\Phi^i\cdot\nabla\Phi^j,\Phi^k)=\int\Phi^i\cdot\nabla\Phi^j\Phi^k=-\int\Phi^i\cdot\nabla\Phi^k\Phi^j=-(\Phi^i\cdot\nabla\Phi^k,\Phi^j).
\end{align*}
Then by summing up~$i,j,k$~respectively from~$1$~to~$m$, we obtain
\begin{align*}
\sum^m_{i,j,k=1}e^k_{ij}\tilde{g}^m_i(t)\tilde{g}^m_j(t)\tilde{g}^m_k(t)=0.
\end{align*}
Consequently, we obtain
\begin{align}
\frac{1}{2}\frac{\mathrm{d}}{\mathrm{d}t}\sum_{k=1}^m(\tilde{g}^m_k(t))^2=-\sum_{i,k=1}^mA_i^k(t)\tilde{g}^m_i(t)\tilde{g}^m_k(t)+\sum^m_{k=1}f^k(t)\tilde{g}_k^m(t).
\end{align}

Denote
\begin{equation*}
|\tilde{g}|^2=\sum_{k=1}^{m}(\tilde{g}_{k}^{m}(t))^2, C_m=\max\limits_{1\leq k\leq m}||\Phi^k(x)||^2_{W^{1,\infty}}.
\end{equation*}

By Cauchy inequality and the assumption (1.12), (3.14) can be rewritten as
\begin{align}
\frac{1}{2}\frac{\mathrm{d}}{\mathrm{d}t}|\tilde{g}|^2\leq C |f|^2+C_m|\tilde{g}|^2,
\end{align}
which yields
\begin{align}
\frac{\mathrm{d}}{\mathrm{d}t}\left(e^{-C_mt}|\tilde{g}|^2\right)\leq C |f|^2,
\end{align}
where~$f=(f_{1}^{m}(t),f_{2}^{m}(t),\cdots,f_{m}^{m}(t)).$

As we know,
\begin{align}
f^k(t)=(\lambda(\theta^m)\nabla d^m\odot \nabla d^m,\nabla \Phi^k)=\int\lambda(\theta^m)\nabla d^m\odot \nabla d^m:\nabla \Phi^k,
\end{align}
so
\begin{align}
(f^k(t))^2=\left|\int\lambda(\theta^m)\nabla d^m\odot \nabla d^m:\nabla \Phi^k\right|^2\leq \bar{\lambda}^2\left(\int|\nabla d^m|^2|\nabla\Phi^k|\right)^2,
\end{align}
which implies
\begin{align}
|f|^2\leq C_m \bar{\lambda}^2||\nabla d^m||^4_{L^2},
\end{align}
where~$C_m=\max\limits_{1\leq k\leq m}||\Phi^k(x)||^2_{W^{1,\infty}}.$

In order to get the estimate of the upper bound of $|f|^2$, one should estimate $||\nabla d^m||^2_{L^2}$.
Since~$d^m(x,t)$~solves
\begin{align}
d^m_t+v^m \cdot\nabla d^m=\Delta d^m-(|d^m|^2-1)d^m,
\end{align}
multiplying~$\Delta d^m-(|d^m|^2-1)d^m$ on the both sides and integrating over $\Omega$, we obtain
\begin{align}
&\frac{1}{2}\frac{\mathrm{d}}{\mathrm{d}t}\int\left(|\nabla d^m|^2+\frac{W(d^m)}{2}\right)+\int|\Delta d^m-(|d^m|^2-1)d^m|^2\nonumber\\
=&\int v^m \cdot\nabla d^m(\Delta d^m-(|d^m|^2-1)d^m).
\end{align}
Applying Cauchy inequality, we derive
\begin{align}
\frac{\mathrm{d}}{\mathrm{d}t}\int\left(|\nabla d^m|^2+\frac{W(d^m)}{2}\right)+\int|\Delta d^m-(|d^m|^2-1)d^m|^2\leq\int|v^m \cdot\nabla d^m|^2.
\end{align}

Noting that
\begin{align}
\|v^m\|^2_{L^{\infty}}&=\|\sum_{i=1}^{m}g_{k}^{m}(t)\Phi^k(x)\|^2_{L^{\infty}}\nonumber\\
&=\sum_{k=1}^{m}|g_{k}^{m}(t)|^2\|\Phi^k(x)\|^2_{L^{\infty}}\nonumber\\
&\leq \sum_{k=1}^{m}|g_{k}^{m}(t)|^2\max\limits_{1\leq k\leq m}\|\Phi^k(x)\|^2_{L^{\infty}}\nonumber\\
&\leq C_m K,
\end{align}
we have
\begin{align}
\frac{\mathrm{d}}{\mathrm{d}t}\int\left(|\nabla d^m|^2+\frac{W(d^m)}{2}\right)+\int|\Delta d^m+(1-|d^m|^2)d^m|^2\leq C_mK\int|\nabla d^m|^2.
\end{align}
Hence, we have reached the following result by applying the Gronwall inequality
\begin{align}
&\sup\limits_{0\leq t \leq T_0}\int\left(|\nabla d^m|^2+\frac{W(d^m)}{2}\right)+\int^{T_0}_{0}\int|\Delta d^m+(1-|d^m|^2)d^m|^2\nonumber\\
\leq &(1+C_m K T_0)e^{C_m K T_0}\int|\nabla d_0|^2.
\end{align}
As a consequence, we obtain $|f|^2\leq C(m,K)$, for any $t \in[0,T_0]$, where~$C(m,K)$ is independent of $T_0$ if we take~$0<T_0\leq 1$~in (3.25).

Now integrating both sides of (3.16) over $[0,t]$, one has
\begin{align}
|\tilde{g}|^2(t)\leq \int^t_0 e^{Ct}(|f|^2+|\tilde{g}(0)|^2)\leq e^{C_mT_0}(C(m,K)T_0+|\tilde{g}(0)|^2),\forall t\in[0,T_0].
\end{align}
Then, taking~$K=e^{C_m}(1+|\tilde{g}(0)|^2)$~and letting~$T_0$ be small enough, we have $C(m,K)T_0\leq 1$, which means that
\begin{align}
|\tilde{g}|^2(t)\leq e^{C_m}(1+|\tilde{g}(0)|^2)=K,\forall t\in[0,T_0].
\end{align}
So, when we take~$T=T_0$, $\widetilde{\mathcal{T}}$ is a continuous map from~$\mathcal{C}_K$ to $\mathcal{C}_K$.

Moreover, the theory of ordinary equation yields that~$\tilde {g}(t)\in(C^1[0,T_0])^m$, which together with the estimation of~$|f|^2$~implies that
\begin{align}
\sum_{k=1}^{m}\left|\frac{d}{dt}\tilde{g}_{k}^{m}(t)\right|^2
&\leq \sum_{k=1}^{m}\left|\sum^m_{i=1}-A_i^k(t)\tilde{g}^m_i(t)\tilde{g}^m_k(t)-\sum^m_{i,j=1}e^k_{ij}\tilde{g}^m_i(t)\tilde{g}^m_j(t)+ f^k(t)\right|^2\nonumber\\
&\leq C\sum_{k=1}^{m}\left[\left|\sum^m_{i=1}A^k_i(t)\tilde{g}^m_i(t)\tilde{g}_k^m(t)\right|^2+\left|\sum_{i,j=1}^{m}e_{ij}^k\tilde{g}_i^m(t)\tilde{g}_j^m(t)\right|^2+\left|f^k(t)\right|^2)\right]\nonumber\\
&\leq C(m,K).
\end{align}

Therefore, by differential mean value theorem and Arzela-Ascoli theorem we finish the proof that $\widetilde{\mathcal{T}}$ is a compact operator. Consequently, according to Schauder fixed point theorem, there exist $\tilde{g}(t)\in\mathcal{C}_K$ such that $\widetilde{\mathcal{T}}(\tilde{g}(t))=\tilde{g}(t)$ is valid.

Reviewing the process of solving this approximate system and, starting from $t=T_0$ again, we conclude that for fixed $m$ and for any time $T>0$, there exists a solution $(u^m(x,t),d^m(x,t), \theta^m(x,t))$ on $Q_{T}=\Omega\times[0,T]$ with $u^m(x,t)\in C^1(Q_{T}), d^m(x,t)\in L^\infty(0,T;H^3), \theta^m(x,t)\in L^\infty(0,T;H^2)$. In fact, as long as the initial data $(u_0,d_0,\theta_0)$ smooth enough, the solution $(u^m,d^m,\theta^m)$ is smooth on $Q_{T}$.

\section{Uniform {\it a priori} estimation}
To show the solvability of the original system (1.8), we should take $m\rightarrow\infty$ in system (3.1). The following estimates are needed.

Before going on, we should point out that by applying Helmholtz decomposition, the equation for the velocity field  can be rewritten as
\begin{align}
u^m_t+u^m\cdot\nabla u^m+\nabla p^m=\mathrm{div}\left[\mu(\theta^m)(\nabla u^m+\nabla^Tu^m)-\lambda(\theta^m)\nabla d^m\odot \nabla d^m\right],
\end{align}
where $\nabla p^m$ is decided by the property that  $\Phi^i,i=1,2,\cdots,m$ are all divergence free. Then we can restate system that $(u^m,d^m,\theta^m)$ solves as follow:
\begin{align}
\begin{cases}
u^m_t+u^m\cdot\nabla u^m+\nabla p^m=\mathrm{div}(\mathbb{S}^m+\sigma^{nd,m})\\
\mathrm{div}u^m=0\\
d^m_t+u^m \cdot\nabla d^m=\Delta d^m+(1-|d^m|^2)d^m\\
\theta^m_t+u^m\cdot\nabla \theta^m=\Delta \theta^m+(\mathbb{S}^m+\sigma^{nd,m}):\nabla u^m.\\
\end{cases}
\end{align}
Here $\mathbb{S}^m=\mu(\theta^m)(\nabla u^m+\nabla^Tu^m),\sigma^{nd,m}=-\lambda(\theta^m)\nabla d^m\odot \nabla d^m$. In this section, we will omit the index "$m$" for simplicity. In detail, We need lemmas as follows.
\begin{lem}
{\rm( F. Lin, C. Wang \cite{16})}For $0<T<\infty$, suppose $u\in L^2(0,T;W),d_0\in H^1(\Omega,\mathbb{R})$ satisfies $|d_0(x)|\leq 1,$ a.e. $\Omega$. If, in addition, $d\in L^2(0,T;H^1(\Omega,\mathbb{R}))$ solves equation $(4.2)_3$, then $|d(x,t)|\leq 1$, a.e. $(x,t)\in \Omega\times[0,T]$.
\end{lem}
\proof For any $k>1$, we define
\begin{equation*}
V_k(x,t)=
\begin{cases}
k^2-1,&\quad |d|>k;\\
|d|^2-1,&\quad 1<|d|\leq k;\\
0,&\quad|d|\leq 1,
\end{cases}
\end{equation*}
then $V_k(x,t)\geq 0$, a.e. $(x,t)\in \Omega\times[0,T]$.
And by equation $(4.2)_3$, we have
\begin{align}
\frac{\mathrm{d}}{\mathrm{d}t}(|d|^2-1)+u\cdot\nabla(|d|^2-1)=\Delta(|d|^2-1)-2\left(|\nabla d|^2+(|d|^2-1)|d|^2\right)
\end{align}
Hence
\begin{align}
\partial_t V_k+u\cdot V_k=\Delta V_k-2\chi_{\{1<|d|\leq k\}}\left(|\nabla d|^2+(|d|^2-1)|d|^2\right),
\end{align}
which further implies
\begin{align}\partial_t V_k+u\cdot V_k\leq \Delta V_k.\end{align}
Multiplying both sides of(4.5) by $V_k$, integrating by parts over $\Omega$, we have
\begin{align}\frac{\mathrm{d}}{\mathrm{d}t}\|V_k\|_{L^2}+\|\nabla V_k\|_{L^2}\leq 0.\end{align}
Then for any $0\leq t\leq T$, integrating on $[0,t]$, we obtain
\begin{align}\|V_k(t)\|^2_{L^2}+\int^t_0\|\nabla V_k(t)\|^2_{L^2}\leq \|V_k(0)\|^2_{L^2}=0,\end{align}
which directly implies that $V_k(x,t)=0$ a.e. $(x,t)\in \Omega\times[0,T]$. Therefore, by the definition of $V_k(x,t)$ we can conclude that
$|d(x,t)|\leq 1$, a.e. $(x,t)\in \Omega\times[0,T]$.
\endproof

\begin{lem}
{\rm(Basic energy inequality)} Let $\Lambda(\cdot)$ solve the Cauchy problem:
\begin{align*}
\begin{cases}
\Lambda^\prime(\cdot)=\lambda(\cdot)^{-1}\\
\Lambda(\underline{\theta}_0)=0.
\end{cases}
\end{align*}
Then, for any time $0\leq t\leq T$, the following inequality is true
\begin{align}
&\sup_{0\leq t\leq T}\int\left(\frac{K+1}{2}|u|^2+(K+1)\theta-\Lambda(\theta)+\frac{W(d)}{4}+\frac{1}{2}|\nabla d|^2\right)\nonumber\\
+&\int^T_0\int\left(\frac{\lambda^\prime(\theta)}{\lambda^2(\theta)}|\nabla\theta|^2+\frac{\mu(\theta)}{2\lambda(\theta)}|\nabla u+\nabla^Tu|^2+|\Delta d-(|d|^2-1)d|^2\right)\nonumber\\
\leq&\int\left(\frac{K+1}{2}|u_0|^2+(K+1)\theta_0-\Lambda(\theta_0)+\frac{1}{2}|\nabla d_0|^2\right),
\end{align}
where $K>0$ is a large enough constant, such that $K\theta-\Lambda(\theta)>0,\ \text{for any}\ \theta\geq \underline{\theta}_0$.
\end{lem}
\proof Firstly, multiplying $(4.2)_1$ by $u$, integrating by parts over $\Omega$ and using $(4.2)_2$, we have
\begin{align}
\frac{1}{2}\frac{\mathrm{d}}{\mathrm{d}t}\int |u|^2=-\int \mu(\theta)(\nabla u+\nabla^Tu):\nabla u+\int\lambda(\theta)\nabla d\odot \nabla d:\nabla u.
\end{align}
Secondly, multiplying $(4.2)_3$ by $\Delta d+(1-|d|^2)d$ and integrating by parts over $\Omega$, we have
\begin{align}
\frac{1}{2}\frac{\mathrm{d}}{\mathrm{d}t}\int\left(|\nabla d|^2+\frac{W(d)}{2}\right)+\int|\Delta d+(1-|d|^2 d)|^2=\int u\cdot\nabla d\cdot \Delta d.
\end{align}
Thirdly, integrating $(4.2)_4$ by parts over $\Omega$, one gets
\begin{align}
\frac{1}{2}\frac{\mathrm{d}}{\mathrm{d}t}\int\theta=\int\mu(\theta)(\nabla u+\nabla^Tu):\nabla u-\int\lambda(\theta)\nabla d\odot\nabla d:\nabla u.
\end{align}
Finally, multiplying $(4.2)_4$ by $\frac{1}{\lambda(\theta)}$ and integrating by parts over $\Omega$, we obtain
\begin{align}
\frac{\mathrm{d}}{\mathrm{d}t}\int\Lambda(\theta)=\int\frac{\lambda^\prime(\theta)}{\lambda^2(\theta)}|\nabla\theta|^2+\int\frac{\mu(\theta)}{2\lambda(\theta)}|\nabla u+\nabla^Tu|^2-\int\nabla d\odot\nabla d:\nabla u.
\end{align}
Then  $(K+1)((4.9)+(4.11))-(4.12)+(4.10)$ gives
\begin{align}
&\frac{\mathrm{d}}{\mathrm{d}t}\int\left(\frac{K+1}{2}|u|^2+(K+1)\theta-\Lambda(\theta)+\frac{W(d)}{4}+\frac{1}{2}|\nabla d|^2\right)\nonumber\\
+&\int\left(\frac{\lambda^\prime(\theta)}{\lambda^2(\theta)}|\nabla\theta|^2+\frac{\mu(\theta)}{2\lambda(\theta)}|\nabla u+\nabla^Tu|^2+|\Delta d-(|d|^2-1)d|^2\right)=0.
\end{align}
Further, integrating for variable $t$ leads to the following estimate
\begin{align}
&\sup_{0\leq t\leq T}\int\left(\frac{K+1}{2}|u|^2+(K+1)\theta-\Lambda(\theta)+\frac{W(d)}{4}+\frac{1}{2}|\nabla d|^2\right)\nonumber\\
+&\int^T_0\int\left(\frac{\lambda^\prime(\theta)}{\lambda^2(\theta)}|\nabla\theta|^2+\frac{\mu(\theta)}{2\lambda(\theta)}|\nabla u+\nabla^Tu|^2+|\Delta d-(|d|^2-1)d|^2\right)\nonumber\\
\leq&\int\left(\frac{K+1}{2}|u_0|^2+(K+1)\theta_0-\Lambda(\theta_0)+\frac{1}{2}|\nabla d_0|^2\right)
\end{align}
where $K$ is a large enough constant, such that $K\theta-\Lambda(\theta)\geq 0$ is valid for any $\theta>\underline{\theta}_0$.
In view of $\lambda^\prime\geq 0,\lambda(\theta)\geq\lambda(\underline{\theta}_0)>0$, we see that $K\geq\lambda(\theta)^{-1}$ is true as long as we take $K=\lambda(\underline{\theta}_0)^{-1}$. Lastly, we integrate it over $[\underline{\theta}_0,\theta]$ to complete the proof of this lemma.
\endproof

Above lemma gives
\begin{align}
&\sup_{0\leq t \leq T}\int\left(|u|^2+2\theta+\frac{W(d)}{2}+|\nabla d|^2\right)+\int^T_0\int\left(\frac{\underline{\mu}}{\bar{\lambda}}|\nabla u|^2+|\Delta d|^2\right)\nonumber\\
\lesssim&\int_{\Omega}\left(|u_0|^2+2(K+1)\theta_0-2\Lambda(\theta_0)+|\nabla d_0|^2\right).
\end{align}
More simply, it yields
\begin{align}
\sup_{0\leq t \leq T}\int(|u|^2+\theta+W(d)+|\nabla d|^2)+\int^T_0\int(|\nabla u|^2+|\Delta d|^2)\leq C,
\end{align}
where $C$ depends only on initial data and known constants.

Moreover, some higher order estimates are proved in the following lemmas.
\begin{lem}
For any $0<t\leq T,\varepsilon>0$, the following inequality is valid.
\begin{align}
&\frac{\mathrm{d}}{\mathrm{d}t}(||\nabla u||^2_{L^2}+||\nabla \theta||^2_{L^2}+||\Delta d||^2_{L^2})+(\underline{\mu}||\Delta u||^2_{L^2}+||\nabla\Delta d||^2_{L^2}+||\Delta \theta||^2_{L^2})\nonumber\\
\leq& C(||\Delta d||^2_{L^2}+||\nabla u||^2_{L^2}+1)(||\Delta u||^2_{L^2}+||\Delta \theta||^2_{L^2}+1)+\varepsilon(||\nabla u||^2_{L^\infty}+||\nabla \theta||^2_{L^\infty}),
\end{align}
where $C$ depends only on initial data and known constants.
\end{lem}
\proof Denote $\mathcal{A}(t)=||\nabla u||^2_{L^2}+||\nabla \theta||^2_{L^2}+||\Delta d||^2_{L^2}$, then
\begin{align*}
\frac{1}{2}\frac{\mathrm{d}}{\mathrm{d}t}\mathcal{A}(t)&=(\nabla u,\nabla u_t)+(\nabla\theta,\nabla\theta_t)+(\Delta d,\Delta d_t)\\
&=-(\Delta u,u_t)-(\Delta \theta,\theta_t)-(\nabla\Delta d,\nabla d_t)\\
&=-(\Delta u,-u\cdot\nabla u+\mathrm{div}(\mu(\theta)(\nabla u+\nabla^Tu)-\lambda(\theta)\nabla d\odot\nabla d))\\
&\quad-(\Delta\theta,-u\cdot\nabla \theta+\Delta\theta+\mu(\theta)(\nabla u+\nabla^Tu):\nabla u-\lambda(\theta)\nabla d\odot\nabla d:\nabla u)\\
&\quad-(\nabla\Delta d,-\nabla(u\cdot\nabla d)+\nabla \Delta d-\nabla((|d|^2-1)d))\\
&=-(\Delta u,\mu(\theta)\Delta u)-(\nabla\Delta d,\nabla\Delta d)-(\Delta\theta,\Delta\theta)\\
&\quad+(\Delta u,u\cdot\nabla u)-(\Delta u,\mu^\prime(\theta)\nabla\theta\cdot(\nabla u+\nabla^Tu))+(\Delta u,\mathrm{div}(\lambda(\theta)\nabla d\odot\nabla d))\\
&\quad+(\Delta\theta,u\cdot\nabla\theta)-(\Delta\theta,\mu(\theta)(\nabla u+\nabla^T u):\nabla u)+(\Delta\theta,\lambda(\theta)\nabla d\odot\nabla d:\nabla u)\\
&\quad+(\nabla\Delta d,\nabla(u\cdot\nabla d))+(\nabla\Delta d,\nabla((|d|^2-1)d))
\end{align*}
Hence, we have
\begin{align}
&\frac{1}{2}\frac{\mathrm{d}}{\mathrm{d}t}\mathcal{A}(t)+\underline{\mu}||\Delta u||^2_{L^2}+||\nabla\Delta d||^2_{L^2}+||\Delta\theta||^2_{L^2}\nonumber\\
\leq&(\Delta u,u\cdot\nabla u)+(\Delta\theta,u\cdot\nabla\theta)+(\nabla\Delta d,\nabla(u\cdot\nabla d))+(\nabla\Delta d,\nabla((|d|^2-1)d))\nonumber\\
&-(\Delta u,\mu^\prime(\theta)\nabla\theta\cdot(\nabla u+\nabla^Tu))+(\Delta u,\mathrm{div}(\lambda(\theta)\nabla d\odot\nabla d))\nonumber\\
&-(\Delta\theta,\mu(\theta)(\nabla u+\nabla^T u):\nabla u)+(\Delta\theta,\lambda(\theta)\nabla d\odot\nabla d:\nabla u)\nonumber\\
=&\sum^8_{i=1}I_i
\end{align}

Next we will apply H\"older inequality, Sobolev embedding inequality, interpolation inequality, Young inequality as well as the fact produced by the periodic boundary conditions, like $||\nabla^2u||^2_{L^2}\leq||\Delta u||^2_{L^2}$, to estimate the above eight terms one by one.

 In the following, when we use the interpolation inequality, we only use the simplest one for simplicity. For example, when we should use the interpolation inequality like $||\nabla u||_{L^4}\lesssim||\nabla u||^{1/4}_{L^2}||\nabla u||^{3/4}_{H^1}$, we simply write $||\nabla u||_{L^4}\lesssim||\nabla u||^{1/4}_{L^2}||\nabla^2 u||^{3/4}_{L^2}$ instead, for the lower order terms brings no difficulty on the proof.
 Now we continue the estimate.
\begin{align}
I_1&\leq||\Delta u||_{L^2}||u||_{L^2}||\nabla u||_{L^\infty}\leq C||\Delta u||^2_{L^2}+\varepsilon||\nabla u||^2_{L^\infty}.\\
I_2&\leq||\Delta \theta||_{L^2}||u||_{L^2}||\nabla \theta||_{L^\infty}\leq C||\Delta \theta||^2_{L^2}+\varepsilon||\nabla \theta||^2_{L^\infty}.\\
I_3&\leq||\nabla\Delta d||_{L^2}||\nabla u||_{L^4}||\nabla d||_{L^4}+||\nabla\Delta d||_{L^2}||\nabla^2 d||_{L^2}||u||_{L^\infty}\nonumber\\
&\lesssim||\nabla\Delta d||_{L^2}||\Delta u||_{L^2}||\Delta d||_{L^2}\nonumber\\
&\leq\varepsilon||\nabla\Delta d||^2_{L^2}+C||\Delta d||^2_{L^2}||\Delta u||^2_{L^2}.\\
I_4&\lesssim||\nabla\Delta d||_{L^2}||\nabla d||_{L^2}\leq \varepsilon||\nabla\Delta d||^2_{L^2}+C.\\
I_5&\lesssim||\nabla\theta||_{L^\infty}||\nabla u||_{L^2}||\Delta u||_{L^2}\leq\varepsilon||\nabla\theta||^2_{L^\infty}+C||\nabla u||^2_{L^2}||\Delta u||^2_{L^2}.\\
I_6&\lesssim||\Delta u||_{L^2}||\nabla\theta||_{L^6}||\nabla d||^2_{L^6}+||\Delta u||_{L^2}||\Delta d||_{L^2}||\nabla d||_{L^\infty}\nonumber\\
&\lesssim||\Delta u||_{L^2}||\Delta \theta||_{L^2}||\Delta d||^2_{L^2}+||\Delta u||_{L^2}||\Delta d||_{L^2}||\nabla\Delta d||_{L^2}\nonumber\\
&\leq\varepsilon||\nabla\Delta d||^2_{L^2}+C(||\Delta u||^2_{L^2}+||\Delta \theta||^2_{L^2})||\Delta d||^2_{L^2}.\\
I_7&\lesssim||\Delta \theta||_{L^2}||\nabla u||_{L^\infty}||\nabla u||_{L^2}\leq \varepsilon||\nabla u||^2_{L^\infty}+C||\Delta \theta||^2_{L^2}||\nabla u||^2_{L^2}.\\
I_8&\lesssim||\Delta \theta||_{L^2}||\nabla u||_{L^6}||\nabla d||^2_{L^6}\leq C||\Delta d||^2_{L^2}(||\Delta u||^2_{L^2}+||\Delta \theta||^2_{L^2}).
\end{align}

Substituting $(4.19)-(4.26)$ into $(4.18)$, we complete the proof of this lemma.
\endproof

\begin{lem}
For any $0\leq t\leq T $, the following estimate is true.
\begin{align}
&\frac{\mathrm{d}}{\mathrm{d}t}(||\Delta u||^2_{L^2}+||\nabla\Delta d||^2_{L^2}+||\Delta \theta||^2_{L^2})+(\underline{\mu}||\nabla\Delta u||^2_{L^2}+||\Delta d_t||^2_{L^2}+||\nabla\Delta \theta||^2_{L^2})\nonumber\\
\leq &C(||\nabla\Delta d||^4_{L^2}+||\Delta u||^4_{L^2}+||\Delta u||^2_{L^2})(||\Delta \theta||^4_{L^2}+||\Delta \theta||^2_{L^2}+||\nabla\Delta d||^2_{L^2}+1)\nonumber\\
&+C||\nabla\Delta d||^4_{L^2}(||\Delta u||^2_{L^2}||\Delta\theta||^2_{L^2}+||\Delta u||^2_{L^2}+||\Delta\theta||^2_{L^2})+C||\Delta d||^2_{L^2},
\end{align}
where $C $ depends only on initial data and known constants.
\end{lem}
\proof Applying $\nabla $ to $(4.2)_{1} $, multiplying the result by $\nabla\Delta u $ and integrating by parts over $\Omega $, we have
\begin{align}
\frac{1}{2}\frac{\mathrm{d}}{\mathrm{d}t}\int|\Delta u|^2=&-\int\nabla\mathrm{div}(\mu(\theta)(\nabla u+\nabla^Tu)):\nabla\Delta u+\int\nabla(u\cdot\nabla u):\nabla\Delta u\nonumber\\
&+\int\nabla[\mathrm{div}(\lambda(\theta)\nabla d\odot\nabla d)]:\nabla\Delta u
\end{align}
Note that
\begin{align*}
\nabla\mathrm{div}(\mu(\theta)(\nabla u+\nabla^Tu))=\mu(\theta)\nabla\Delta u+\mu^\prime(\theta)\nabla\theta\Delta u-\nabla[\mu^\prime(\theta)\nabla\theta\cdot(\nabla u+\nabla^Tu)],
\end{align*}
hence (4.28) can be rewritten as
\begin{align}
&\frac{1}{2}\frac{\mathrm{d}}{\mathrm{d}t}\int|\Delta u|^2+\int\mu(\theta)|\nabla\Delta u|^2\nonumber\\
=&\int\mu^\prime(\theta)\nabla\theta\otimes\Delta u:\nabla\Delta u-\int\nabla[\mu^\prime(\theta)\nabla\theta\cdot(\nabla u+\nabla^Tu)]:\nabla\Delta u\nonumber\\
&+\int\nabla(u\cdot\nabla u):\nabla\Delta u+\int\nabla[\mathrm{div}(\lambda(\theta)\nabla d\odot\nabla d)]:\nabla\Delta u.
\end{align}
Next, applying $\Delta $ to $(4.2)_3$, multiplying the result by $\Delta d_{t} $ and integrating by parts over $\Omega $, we obtain
\begin{align}
\frac{1}{2}\frac{\mathrm{d}}{\mathrm{d}t}\int|\nabla\Delta d|^2+\int|\Delta d_t|^2=-\int\Delta[(|d|^{2}-1)d]\cdot\Delta d_t-\int\Delta(u\cdot\nabla d)\cdot\Delta d_t,
 \end{align}
while applying $\nabla $ to $(4.2)_{4} $, multiplying the result by $\nabla\Delta \theta $ and integrating by parts over $\Omega $, one gets
\begin{align}
\frac{1}{2}\frac{\mathrm{d}}{\mathrm{d}t}\int|\Delta \theta|^2+\int|\nabla\Delta \theta|^2=&\int\nabla(u\cdot\nabla \theta):\nabla\Delta \theta+\int\nabla[\lambda(\theta)\nabla d\odot\nabla d:\nabla u]\cdot\nabla\Delta\theta\nonumber\\
&-\int\nabla(\mu(\theta)(\nabla u+\nabla^Tu):\nabla u)\cdot\nabla\Delta\theta.
\end{align}
Putting (4.29), (4.30) and (4.31) together, we have
\begin{align}
&\frac{1}{2}\frac{\mathrm{d}}{\mathrm{d}t}(||\Delta u||^2_{L^2}+||\nabla\Delta d||^2_{L^2}+||\Delta \theta||^2_{L^2})+(\underline{\mu}||\nabla\Delta u||^2_{L^2}+
||\Delta d_t||^2_{L^2}+||\nabla\Delta \theta||^2_{L^2}) \nonumber\\
\leq&\int\mu^\prime(\theta)\nabla\theta\otimes\Delta u:\nabla\Delta u-\int\nabla[\mu^\prime(\theta)\nabla\theta\cdot(\nabla u+\nabla^Tu)]:\nabla\Delta u\nonumber\\
&+\int\nabla(u\cdot\nabla u):\nabla\Delta u+\int\nabla[\mathrm{div}(\lambda(\theta)\nabla d\odot\nabla d)]:\nabla\Delta u\nonumber\\
&-\int\Delta[(|d|^{2}-1)d]\cdot\Delta d_t-\int\Delta(u\cdot\nabla d)\cdot\Delta d_t+\int\nabla(u\cdot\nabla \theta):\nabla\Delta \theta\nonumber\\
&+\int\nabla[\lambda(\theta)\nabla d\odot\nabla d:\nabla u]\cdot\nabla\Delta\theta-\int\nabla(\mu(\theta)(\nabla u+\nabla^Tu):\nabla u)\cdot\nabla\Delta\theta\nonumber\\
=&\sum^9_{i=1}J_i.
\end{align}

Similarly, when we use interpolation inequality or Sobolev embedding inequality, we only calculate the highest order terms for convenience. We estimate the above nine terms one by one as follows:
\begin{align}
J_1&\lesssim\int|\nabla\theta||\Delta u||\nabla\Delta u|\lesssim||\nabla\Delta u||_{L^2}||\nabla\theta||_{L^6}||\Delta u||_{L^3}\nonumber\\
&\lesssim||\nabla\Delta u||^{\frac{3}{2}}_{L^2}||\Delta\theta||_{L^2}||\Delta u||^{\frac{1}{2}}_{L^2}\leq\varepsilon||\nabla\Delta u||^2_{L^2}+C||\Delta\theta||^4_{L^2}||\Delta u||^2_{L^2}.\\
J_2&=\int\mu^\prime(\theta)\partial_k\partial_i\theta(\partial_iu^j+\partial_ju^i)\partial_k\Delta u^j
+\int\mu^\prime(\theta)\partial_i\theta(\partial_k\partial_iu^j+\partial_k\partial_ju^i)\partial_k\Delta u^j\nonumber\\
&\quad+\int\mu^{\prime\prime}(\theta)\partial_k\theta\partial_i\theta(\partial_iu^j+\partial_ju^i)\partial_k\Delta u^j\nonumber\\
&\lesssim\int|\nabla^2\theta||\nabla u||\nabla\Delta u|+\int|\nabla\theta||\nabla^2u||\nabla\Delta u|+\int|\nabla\theta|^2|\nabla u||\nabla\Delta u|\nonumber\\
&\lesssim||\nabla\Delta u||_{L^2}(||\nabla\theta||_{L^6}||\nabla^2u||_{L^3}+||\nabla u||_{L^6}||\nabla^2\theta||_{L^3}+||\nabla u||_{L^6}||\nabla\theta||^2_{L^6})\nonumber\\
&\lesssim||\nabla\Delta u||^\frac{3}{2}_{L^2}||\Delta u||^\frac{1}{2}_{L^2}||\Delta\theta||_{L^2}+||\nabla\Delta u||_{L^2}||\Delta u||_{L^2}||\nabla\Delta\theta||^\frac{1}{2}_{L^2}||\Delta\theta||^\frac{1}{2}_{L^2}\nonumber\\
&\quad+||\nabla\Delta u||_{L^2}||\Delta u||_{L^2}||\Delta\theta||^2_{L^2}\nonumber\\
&\leq\varepsilon||\nabla\Delta u||^2_{L^2}+\varepsilon||\nabla\Delta\theta||^2_{L^2}+C(||\Delta\theta||^2_{L^2}+||\Delta\theta||^4_{L^2})(||\Delta u||^2_{L^2}+||\Delta u||^4_{L^2}).\\
J_3&=\int\partial_i(u^k\partial_k u^j)\partial_i\Delta u^j=\int\partial_i u^k\partial_k u^j\partial_i\Delta u^j+\int u^k\partial_i\partial_k u^j\partial_i\Delta u^j\nonumber\\
&\lesssim\int(|\nabla u|^2|\nabla\Delta u|+|\nabla^2 u||u||\nabla\Delta u|)\nonumber\\
&\lesssim||\nabla\Delta u||_{L^2}||\nabla u||^2_{L^4}+||\nabla^2 u||_{L^2}|| u||_{L^\infty}\nonumber\\
&\leq \varepsilon||\nabla\Delta u||^2_{L^2}+C||\Delta u||^4_{L^2}.\\
J_4&=\int\partial_l\partial_i(\lambda(\theta)\partial_id\cdot\partial_jd)\partial_l\Delta u^j\nonumber\\
&=\int\partial_l[\lambda^\prime(\theta)\partial_i\theta\partial_i d\cdot\partial_j d+\lambda(\theta)\Delta d\cdot\partial_j d+\lambda(\theta)\partial_i d\cdot\partial_j\partial_id]
\partial_l\Delta u^j\nonumber\\
&=\int[\lambda^{\prime\prime}(\theta)\partial_l\theta\partial_i\theta\partial_id\cdot\partial_jd+\lambda^\prime(\theta)\partial_l\partial_i\theta\partial_id\cdot\partial_jd
+\lambda^\prime(\theta)\partial_i\theta\partial_l(\partial_id\cdot\partial_jd)]\partial_l\Delta^j\nonumber\\
&\quad+\int[\lambda^\prime(\theta)\partial_l\theta\Delta d\cdot\partial_jd+\lambda(\theta)\partial_l\Delta d\cdot\partial_jd+\lambda(\theta)\Delta d\cdot\partial_l\partial_jd]
\partial_l\Delta u^j\nonumber\\
&\quad+\int[\lambda^\prime(\theta)\partial_l\theta\partial_id\cdot\partial_i\partial_jd+\lambda(\theta)\partial_l(\partial_id\cdot\partial_i\partial_jd)]\partial_l\Delta u^j.\nonumber
\end{align}
Noting that
\begin{align*}
&\int\lambda(\theta)\Delta d\cdot\partial_l\partial_jd\partial_l\Delta u^j=\int-\lambda^\prime(\theta)\partial_j\theta\Delta d\cdot\partial_ld\partial_l\Delta u^j
-\lambda(\theta)\partial_j\Delta d\cdot\partial_l d\partial_l\Delta u^j,\\
&\int\lambda(\theta)\partial_l(\partial_id\cdot\partial_i\partial_jd)\partial_l\Delta u^j=-\frac{1}{2}\int\lambda^\prime(\theta)\partial_j\theta\partial_l(|\nabla d|^2)\partial_l\Delta u^j,
\end{align*}
we have
\begin{align}
J_4&\lesssim\int(|\nabla\theta|^2|\nabla d|^2+|\nabla^2\theta||\nabla d|^2+|\nabla\theta||\nabla^2 d||\nabla d|+|\nabla\Delta d||\nabla d|)|\nabla\Delta u|\nonumber\\
&\lesssim||\nabla\Delta u||_{L^2}||\nabla d||^2_{L^\infty}(||\nabla\theta||^2_{L^4}+||\nabla^2\theta||_{L^2})\nonumber\\
&\quad+||\nabla\Delta u||^2_{L^2}||\nabla d||_{L^\infty}(||\nabla\theta||_{L^3}||\nabla^2 d||_{L^6}+||\nabla\Delta d||_{L^2})\nonumber\\
&\leq \varepsilon||\nabla\Delta u||^2_{L^2}+C||\nabla\Delta d||^4_{L^2}(||\Delta \theta||^4_{L^2}+||\Delta \theta||^2_{L^2}+1).\\
J_5&=-\int\Delta(|d|^2-1)d\cdot\Delta d_t-\int(|d|^2-1)\Delta d\cdot\Delta d_t-\int 2\nabla(|d|^2-1)\nabla d\cdot\Delta d_t\nonumber\\
&\lesssim\int(|\Delta d||\Delta d_t|+|\nabla d|^2|\Delta d_t|)\nonumber\\
&\leq\varepsilon||\Delta d_t||^2_{L^2}+C(||\Delta d||^2_{L^2}+||\nabla d||^4_{L^2})\nonumber\\
&\leq\varepsilon||\Delta d_t||^2_{L^2}+C||\Delta d||^2_{L^2}.\\
J_6&=-\int(\Delta u\cdot\nabla d\cdot\Delta d_t+u\cdot\nabla\Delta d\cdot\Delta d_t+2\nabla u\nabla^2d\cdot\Delta d_t)\nonumber\\
&\lesssim\int(|\Delta u||\nabla d|+|u||\nabla\Delta d|+|\nabla u||\nabla^2d|)|\Delta d_t|\nonumber\\
&\lesssim||\Delta d_t||_{L^2}(||\Delta u||_{L^2}||\nabla d||_{L^\infty}+||\nabla\Delta d||_{L^2}||u||_{L^\infty}+||\nabla u||_{L^3}||\nabla^2 d||_{L^6})\nonumber\\
&\leq\varepsilon||\Delta d_t||^2_{L^2}+C||\Delta u||^2_{L^2}||\nabla\Delta d||^2_{L^2}.\\
J_7&=\int\partial_i(u^k\partial_k\theta)\partial_i\Delta\theta
=\int(\partial_i u^k\partial_k\theta\partial_i\Delta\theta+u^k\partial_i\partial_k\theta\partial_i\Delta\theta)\nonumber\\
&\leq\varepsilon||\nabla\Delta \theta||^2_{L^2}+C||\nabla u||^2_{L^4}||\nabla\theta||^2_{L^4}+C|| u||^2_{L^\infty}||\Delta \theta||_{L^2}\nonumber\\
&\leq\varepsilon||\nabla\Delta \theta||^2_{L^2}+C||\Delta u||^2_{L^2}||\Delta\theta||^2_{L^2}.\\
J_8&=\int\partial_i[\lambda(\theta)\partial_id\cdot\partial_jd\partial_iu^j]\partial_k\Delta\theta\nonumber\\
&=\int\lambda^\prime(\theta)\partial_k\theta\partial_id\cdot\partial_jd\partial_iu^j\partial_k\Delta\theta
+\int\lambda(\theta)\partial_k\partial_id\cdot\partial_jd\partial_iu^j\partial_k\Delta\theta\nonumber\\
&\quad +\int\lambda(\theta)\partial_id\cdot\partial_k\partial_jd\partial_iu^j\partial_k\Delta\theta
+\int\lambda(\theta)\partial_id\cdot\partial_jd\partial_k\partial_iu^j\partial_k\Delta\theta\nonumber\\
&\lesssim\int|\nabla\Delta\theta|(|\nabla\theta||\nabla d|^2|\nabla u|+|\nabla^2 d||\nabla d||\nabla u|+|\nabla d|^2|\nabla^2 u|)\nonumber\\
&\lesssim (||\nabla d||^2_{L^\infty}||\nabla\theta||_{L^4}||\nabla u||_{L^4}+||\nabla^2 d||_{L^4}||\nabla d||_{L^\infty}||\nabla u||_{L^4}
+||\nabla d||^2_{L^\infty}||\nabla^2 u||_{L^2})||\nabla\Delta \theta||_{L^2}\nonumber\\
&\lesssim||\nabla\Delta \theta||_{L^2}||\nabla\Delta d||^2_{L^2}||\Delta u||_{L^2}(||\Delta \theta||_{L^2}+1)\nonumber\\
&\leq\varepsilon||\nabla\Delta \theta||^2_{L^2}+C||\nabla\Delta d||^4_{L^2}||\Delta u||^2_{L^2}||\Delta \theta||^2_{L^2}+C||\nabla\Delta d||^4_{L^2}||\Delta u||^2_{L^2}.\\
J_9&=\int\mu^\prime(\theta)\partial_k\theta(\partial_iu^j+\partial_ju^i)\partial_iu^j\partial_k\Delta\theta
+\int\mu(\theta)(\partial_k\partial_iu^j+\partial_k\partial_ju^i)\partial_iu^j\partial_k\Delta\theta\nonumber\\
&\quad+\int\mu(\theta)(\partial_iu^j+\partial_ju^i)\partial_k\partial_iu^j\partial_k\Delta\theta\nonumber\\
&\lesssim\int|\nabla\theta||\nabla u|^2|\nabla\Delta\theta|+\int|\nabla u||\nabla^2 u||\nabla\Delta\theta|\nonumber\\
&\lesssim||\nabla\Delta \theta||_{L^2}||\nabla u||^2_{L^6}||\nabla\theta||_{L^6}+||\nabla\Delta \theta||_{L^2}||\nabla^2 u||_{L^3}||\nabla u||_{L^6}\nonumber\\
&\lesssim||\nabla\Delta \theta||_{L^2}||\Delta u||^2_{L^2}||\Delta\theta||_{L^2}+||\nabla\Delta \theta||_{L^2}||\Delta u||^\frac{3}{2}_{L^2}||\nabla\Delta u||^\frac{1}{2}_{L^2}\nonumber\\
&\leq\varepsilon||\nabla\Delta \theta||^2_{L^2}+\varepsilon||\nabla\Delta u||^2_{L^2}+C||\Delta u||^4_{L^2}(||\Delta\theta||^2_{L^2}+||\Delta u||^2_{L^2}).
\end{align}
At last, substituting (4.33)-(4.41) into (4.32), we obtain
\begin{align}
&\frac{1}{2}\frac{\mathrm{d}}{\mathrm{d}t}(||\Delta u||^2_{L^2}+||\nabla\Delta d||^2_{L^2}+||\Delta \theta||^2_{L^2})+(\underline{\mu}||\nabla\Delta u||^2_{L^2}+||\Delta d_t||^2_{L^2}+||\nabla\Delta \theta||^2_{L^2})\nonumber\\
\leq &C_0\varepsilon(||\nabla\Delta u||^2_{L^2}+||\Delta d_t||^2_{L^2}+||\nabla\Delta\theta||^2_{L^2})+C||\Delta d||^2_{L^2}\nonumber\\
&+C(||\nabla\Delta d||^4_{L^2}+||\Delta u||^4_{L^2}+||\Delta u||^2_{L^2})(||\Delta \theta||^4_{L^2}+||\Delta \theta||^2_{L^2}+||\nabla\Delta d||^2_{L^2}+1)\nonumber\\
&+C||\nabla\Delta d||^4_{L^2}(||\Delta u||^2_{L^2}||\Delta\theta||^2_{L^2}+||\Delta u||^2_{L^2}+||\Delta\theta||^2_{L^2}).
\end{align}
We can take $\varepsilon$ small enough such that $C_0\varepsilon\leq\frac{1}{2}\min\{\underline{\mu},1\}$. The lemma follows.
\endproof

Thanks to lemmas above, we obtain the local uniform estimate to the approximate solution. That is
\begin{lem}
There exists a finite time $T^*>0$ such that for any $0<T<T^*$, the following inequality is true.
\begin{align}
&\sup_{0\leq t\leq T}(||\nabla u||^2_{H^1}+||\nabla \theta||^2_{H^1}+||\Delta d||^2_{H^1})\nonumber\\
&+\int^T_0(||\Delta u||^2_{H^1}+||\Delta \theta||^2_{H^1}+||\nabla\Delta d||^2_{L^2}+||\Delta d_t||^2_{L^2})\leq C,
\end{align}
where $C$ depends only on initial data and known constants.
\end{lem}
\proof
Adding (4.17) to (4.27), one has
\begin{align}
&\frac{\mathrm{d}}{\mathrm{d}t}(||\nabla u||^2_{H^1}+||\nabla \theta||^2_{H^1}+||\Delta d||^2_{H^1})+(\underline{\mu}||\Delta u||^2_{H^1}+||\Delta \theta||^2_{H^1}+||\nabla\Delta d||^2_{L^2}+||\Delta d_t||^2_{L^2})\nonumber\\
\leq& C(||\Delta d||^2_{L^2}+||\nabla u||^2_{L^2}+1)(||\Delta u||^2_{L^2}+||\Delta \theta||^2_{L^2}+1)+\varepsilon(||\nabla u||^2_{L^\infty}+||\nabla \theta||^2_{L^\infty})\nonumber\\
&+C(||\nabla\Delta d||^4_{L^2}+||\Delta u||^4_{L^2}+||\Delta u||^2_{L^2})(||\Delta \theta||^4_{L^2}+||\Delta \theta||^2_{L^2}+||\nabla\Delta d||^2_{L^2}+1)\nonumber\\
&+C||\nabla\Delta d||^4_{L^2}(||\Delta u||^2_{L^2}||\Delta\theta||^2_{L^2}+||\Delta u||^2_{L^2}+||\Delta\theta||^2_{L^2}).
\end{align}
By Sobolev embedding inequality, there exists a constant $C>0$, such that
\begin{align*}
||\nabla u||_{L^\infty}+||\nabla \theta||_{L^\infty}\leq C(||\nabla u||_{H^2}+||\nabla \theta||_{H^2}).
\end{align*}
Substituting the result into (4.44) and taking $\varepsilon$ small enough,  we have
\begin{align}
&\frac{\mathrm{d}}{\mathrm{d}t}(||\nabla u||^2_{H^1}+||\nabla \theta||^2_{H^1}+||\Delta d||^2_{H^1})+(||\Delta u||^2_{H^1}+||\Delta \theta||^2_{H^1}+||\nabla\Delta d||^2_{L^2}+||\Delta d_t||^2_{L^2})\nonumber\\
\leq& C(||\Delta d||^2_{L^2}+||\nabla u||^2_{L^2}+1)(||\Delta u||^2_{L^2}+||\Delta \theta||^2_{L^2}+1))\nonumber\\
&+C(||\nabla\Delta d||^4_{L^2}+||\Delta u||^4_{L^2}+||\Delta u||^2_{L^2})(||\Delta \theta||^4_{L^2}+||\Delta \theta||^2_{L^2}+||\nabla\Delta d||^2_{L^2}+1)\nonumber\\
&+C||\nabla\Delta d||^4_{L^2}(||\Delta u||^2_{L^2}||\Delta\theta||^2_{L^2}+||\Delta u||^2_{L^2}+||\Delta\theta||^2_{L^2}).
\end{align}
Denoting
\begin{align*}
\mathcal{F}(t)&:=||\nabla u||^2_{H^1}+||\nabla \theta||^2_{H^1}+||\Delta d||^2_{H^1}+1,\\
\mathcal{H}(t)&:=||\Delta u||^2_{H^1}+||\Delta \theta||^2_{H^1}+||\nabla\Delta d||^2_{L^2}+||\Delta d_t||^2_{L^2},
\end{align*}
one can rewrite (4.41) as
\begin{align}
\frac{\mathrm{d}}{\mathrm{d}t}\mathcal{F}(t)+\mathcal{H}(t)\leq C\mathcal{F}(t)^4
\end{align}

In particular, we have
\begin{align}
\frac{\mathrm{d}}{\mathrm{d}t}\mathcal{F}(t)\leq C\mathcal{F}(t)^4.
\end{align}
Again, we rewrite the above inequality as
\begin{align*}
\frac{\mathrm{d}}{\mathrm{d}t}\mathcal{F}(t)^{-3}\geq -3C.
\end{align*}
Then integrating both sides of the inequality over $[0,t]$, we obtain
\begin{align*}
\mathcal{F}(t)^{-3}\geq-3Ct
\end{align*}
Here,we denote $T^*:=[3C\mathcal{F}(0)^3]^{-1}$. Therefore, for any $T>0$, it derives
\begin{align}
\mathcal{F}(t)^3\leq[\mathcal{F}(0)^3-3CT]^{-1},\forall 0<t\leq T,
\end{align}
as long as $T<T^*$. Finally, for any $t\in(0,T)$, we integrate (4.42) over $[0,t]$ to complete the proof.
\endproof

For a more complete description to the approximate solution, we need the following lemma.
\begin{lem}
For any $0<T<T^*$, it is valid that
\begin{align}
\sup_{0\leq t\leq T}||p||^2_{H^1}+\int^T_0(||p||^2_{H^2}+||\nabla u_t||^2_{L^2}+||\nabla \theta_t||^2_{L^2}+||\Delta d||^2_{H^2})\leq C.
\end{align}
Here $T^*$ is the finite time decided by Lemma 4.5 and $C$ depends only on initial data and known constants.
\end{lem}
\proof
Similar to Lemma 4.3, when we use interpolation inequality or Sobolev embedding inequality, we only calculate the highest order terms for convenience.

First of all, applying $\mathrm{div}$ to $(4.2)_1$, by $\mathrm{div}u=0$ we have
\begin{align}
\Delta p=\mathrm{divdiv}(\mu(\theta(\nabla u+\nabla^Tu))-\lambda(\theta)\nabla d\odot\nabla d)-\mathrm{div}(u\cdot\nabla u).
\end{align}
Thus, by Lemma 4.5 and the standard elliptic estimates, together with H\"older inequality, Sobolev embedding inequality and Young inequality, we have
\begin{align}
||p||^2_{H^1}&\lesssim||\mathrm{div}(\mu(\theta)(\nabla u+\nabla^Tu))||^2_{L^2}+||\mathrm{div}(\lambda(\theta)\nabla d\odot\nabla d)||^2_{L^2}+||u\cdot\nabla u||^2_{L^2}\nonumber\\
&\lesssim\int|\mu^\prime(\theta)\nabla\theta\cdot(\nabla u+\nabla^Tu)|^2+\int|\mu(\theta)\Delta u|^2+\int|u\cdot\nabla u|^2\nonumber\\
&\quad+\int|\lambda^\prime(\theta)\nabla\theta\nabla d\odot\nabla d|^2+\int|\lambda(\theta)\Delta d\cdot\nabla d|^2+\int|\lambda(\theta)\nabla d\:\nabla^2 d|^2\nonumber\\
&\lesssim\int(|\nabla\theta|^2|\nabla u|^2+|\Delta u|^2+|\nabla\theta|^2|\nabla d|^4+|\nabla d|^2|\nabla^2 d|^2+|u|^2|\nabla u|^2)\nonumber\\
&\lesssim||\nabla\theta||^2_{L^4}||\nabla u||^2_{L^4}+||\Delta u||^2_{L^2}+||u||^2_{L^\infty}||\nabla u||^2_{L^2}\nonumber\\
&\quad+||\nabla d||^4_{L^\infty}||\nabla\theta||^2_{L^2}+||\nabla d||^2_{L^\infty}||\Delta d||^2_{L^2}\nonumber\\
&\lesssim||\Delta u||^2_{L^2}||\Delta\theta||^2_{L^2}+||\Delta u||^2_{L^2}+||\Delta u||^2_{L^2}||\nabla u||^2_{L^2}\nonumber\\
&\quad+||\nabla\Delta d||^4_{L^2}||\nabla \theta||^2_{L^2}+||\nabla\Delta d||^2_{L^2}||\Delta d||^2_{L^2},\\
||p||^2_{H^2}&\lesssim||\mathrm{divdiv}(\mu(\theta)(\nabla u+\nabla^Tu))||^2_{L^2}+||\mathrm{divdiv}(\lambda(\theta)\nabla d\cdot\nabla d)||^2_{L^2}+||\mathrm{div}(u\cdot\nabla u)||^2_{L^2}\nonumber\\
&\lesssim\int(|\nabla\theta|^4|\nabla u|^2+|\nabla^2\theta|^2|\nabla u|^2+|\nabla\theta|^2|\Delta u|^2)+\int(|\nabla\theta|^4|\nabla d|^4+|\nabla^2\theta|^2|\nabla d|^4)\nonumber\\
&\quad+\int(|\nabla\theta|^2|\nabla^2d|^2|\nabla d|^2+|\nabla\Delta d|^2|\nabla d|^2+|\nabla^2 d|^4)+\int|\nabla u|^4\nonumber\\
&\lesssim(||\nabla \theta||^4_{L^4}+||\Delta \theta||^2_{L^2})||\nabla u||^2_{L^\infty}+||\nabla \theta||^2_{L^\infty}(||\Delta u||^2_{L^2}+||\nabla d||^4_{L^\infty}||\nabla\theta||^2_{L^2})\nonumber\\
&\quad+||\nabla d||^4_{L^\infty}||\Delta\theta||^2_{L^2}+||\nabla d||^2_{L^\infty}||\nabla \theta||^2_{L^\infty}||\nabla^2 d||^2_{L^2}+||\nabla d||^2_{L^\infty}||\nabla\Delta d||^2_{L^2}\nonumber\\
&\quad+||\nabla^2 d||^4_{L^4}+||\nabla u||^2_{L^2}||\nabla u||^2_{L^\infty}\nonumber\\
&\lesssim(||\Delta \theta||^4_{L^2}+||\Delta \theta||^2_{L^2}+1)||\nabla\Delta u||^2_{L^2}+||\nabla\Delta d||^4_{L^2}||\Delta \theta||^2_{L^2}+||\nabla\Delta d||^4_{L^2}\nonumber\\
&\quad+(||\Delta u||^2_{L^2}+||\nabla\Delta d||^4_{L^2}||\nabla \theta||^2_{L^2}+||\nabla\Delta d||^2_{L^2})||\nabla\Delta \theta||^2_{L^2}\nonumber\\
&\lesssim||\nabla\Delta \theta||^2_{L^2}+||\nabla\Delta u||^2_{L^2}+1.
\end{align}
Therefore, Lemma 4.5 immediately implies that
\begin{align}
\sup_{0\leq t\leq T}||p||^2_{H^1}+\int^T_0||p||^2_{H^2}\leq C.
\end{align}

Secondly, applying $\nabla$ to $(4.2)_1$, we have
\begin{align}
\nabla u_t=\nabla\mathrm{div}(\mu(\theta)(\nabla u+\nabla^Tu))-\nabla\mathrm{div}(\lambda(\theta)\nabla d\odot\nabla d)-\nabla(u\cdot\nabla u)-\nabla^2 p,
\end{align}
which implies that
\begin{align}
||\nabla u_t||^2_{L^2}\lesssim&||\nabla\mathrm{div}(\mu(\theta)(\nabla u+\nabla^Tu))||^2_{L^2}+||\nabla(u\cdot\nabla u)||^2_{L^2}\nonumber\\
&+||\nabla^2 p||^2_{L^2}+||\nabla\mathrm{div}(\lambda(\theta)\nabla d\odot\nabla d)||^2_{L^2},
\end{align}
where
\begin{align}
||\nabla(u\cdot\nabla u)||^2_{L^2}&\lesssim\int|\nabla u|^4+\int|u|^2|\nabla^2 u|^2\nonumber\\
&\lesssim||\nabla u||^2_{L^2}||\nabla u||^2_{L^\infty}+||u||^2_{L^4}||\nabla^2 u)||^2_{L^4}\nonumber\\
&\lesssim||\nabla\Delta u)||^2_{L^2}.
\end{align}
Similarly to $(4.52)$, we have
\begin{align}
&||\nabla\mathrm{div}(\mu(\theta)(\nabla u+\nabla^Tu))||^2_{L^2}+||\nabla\mathrm{div}(\lambda(\theta)\nabla d\odot\nabla d)||^2_{L^2}\nonumber\\
\lesssim &||\nabla\Delta \theta||^2_{L^2}+||\nabla\Delta u||^2_{L^2}+1.
\end{align}
Now, substituting $(4.56),(4.57)$ into $(4.55)$, we obtain
\begin{align}
\int_0^T||\nabla u_t||^2_{L^2}\leq C,
\end{align}
where (4.43) is used.

Thirdly, applying $\nabla$ to $(4.2)_4$ we derive that
\begin{align}
\nabla\theta_t=\nabla\Delta\theta-\nabla(u\cdot\nabla\theta)+\nabla(\mu(\theta)(\nabla u+\nabla^T u):\nabla u)+\nabla(\lambda(\theta)\nabla d\odot\nabla d:\nabla u),
\end{align}
which means
\begin{align}
||\nabla\theta_t||^2_{L^2}\lesssim &||\nabla\Delta\theta||^2_{L^2}+||\nabla(u\cdot\nabla\theta)||^2_{L^2}+||\nabla(\mu(\theta)(\nabla u+\nabla^Tu):\nabla u )||^2_{L^2}\nonumber\\
&+||\nabla(\lambda(\theta)\nabla d\odot\nabla d:\nabla u)||^2_{L^2},
\end{align}
where
\begin{align}
||\nabla(u\cdot\nabla\theta)||^2_{L^2}&\lesssim\int|\nabla u|^2|\nabla\theta|^2+\int|u|^2|\nabla^2\theta|^2\nonumber\\
&\lesssim||\nabla u||^2_{L^2}||\nabla \theta||^2_{L^\infty}+|| u||^2_{L^4}||\nabla^2 \theta||^2_{L^4}\nonumber\\
&\lesssim||\nabla\Delta \theta||^2_{L^2},\\
||\nabla(\mu(\theta)(\nabla u+\nabla^T u):\nabla u)||^2_{L^2}&\lesssim\int|\nabla\theta|^2|\nabla u|^4+\int|\nabla^2u|^2|\nabla u|^2\nonumber\\
&\lesssim||\nabla u||^4_{L^4}||\nabla\theta||^2_{L^\infty}+||\nabla^2 u||^2_{L^2}||\nabla u||^2_{L^\infty}\nonumber\\
&\lesssim||\nabla\Delta u||^2_{L^2}+||\nabla\Delta\theta||^2_{L^2},\\
||\nabla(\lambda(\theta)\nabla d\odot\nabla d:\nabla u)||^2_{L^2}&\lesssim\int(|\nabla\theta|^2|\nabla d|^4|\nabla u|^2+|\nabla^2d|^2|\nabla d|^2|\nabla u|^2+|\nabla d|^4|\nabla^2 u|^2)\nonumber\\
&\lesssim||\nabla d||^2_{L^\infty}||\nabla \theta||^2_{L^\infty}||\nabla d||^2_{L^2}+||\nabla d||^2_{L^\infty}||\nabla u||^2_{L^\infty}||\nabla^2 d||^2_{L^2}\nonumber\\
&\quad+||\nabla d||^4_{L^\infty}||\nabla^2 u||^2_{L^2}\nonumber\\
&\lesssim||\nabla\Delta d||^2_{L^2}(||\nabla\Delta u||^2_{L^2}+||\nabla\Delta \theta||^2_{L^2}+||\nabla\Delta d||^2_{L^2})\nonumber\\
&\lesssim||\nabla\Delta u||^2_{L^2}+||\nabla\Delta \theta||^2_{L^2}+1.
\end{align}
Substituting $(4.61)-(4.63)$ into $(4.60)$ and using (4.43), we have
\begin{align}
\int_0^T||\nabla \theta_t||^2_{L^2}\leq C.
\end{align}
Lastly, applying $\Delta$ to $(4.2)_3$, we derive that
\begin{align}
\Delta^2d=\Delta d_t+\Delta(u\cdot\nabla d)+\Delta[(|d|^2-1)d].
\end{align}
By Lemma 4.5 and the standard elliptic estimates, together with H\"older inequality, Sobolev embedding inequality and Young inequality, we have
\begin{align}
||\Delta d||^2_{H^2}\lesssim||\Delta d_t||^2_{L^2}+||\Delta (u\cdot\nabla d)||^2_{L^2}+||\Delta[(|d|^2-1)d]||^2_{L^2},
\end{align}
where we note that
\begin{align}
||\Delta (u\cdot\nabla d)||^2_{L^2}&\lesssim\int(|u|^2|\nabla\Delta d|^2+|\nabla u|^2|\nabla^2d|^2+|\Delta u|^2|\nabla d|^2)\nonumber\\
&\lesssim||\nabla\Delta d||^2_{L^2}||u||^2_{L^\infty}+||\nabla u||^2_{L^\infty}||\nabla^2 d||^2_{L^2}+||\nabla d||^2_{L^\infty}||\Delta u||^2_{L^2}\nonumber\\
&\lesssim||\nabla\Delta d||^2_{L^2}||\Delta d||^2_{L^2}+||\nabla\Delta u||^2_{L^2}\nonumber\\
&\lesssim||\nabla\Delta u||^2_{L^2}+1,\\
||\Delta[(|d|^2-1)d]||^2_{L^2}&\lesssim\int|\nabla d|^4+\int|\Delta d|^2\leq C.
\end{align}
Then substituting $(4.67)-(4.68)$ into $(4.66)$ and using (4.43), we obtain
\begin{align}
\int_0^T||\Delta d||^2_{H^2}\leq C.
\end{align}
So far, we complete the proof.
\endproof

\section{Existence of the local strong solutions}
Concluding from the previous section, we have shown that there exist $T^*>0$, such that for any $T\in(0,T^*)$, the following estimates are valid.
\begin{align}
&\sup_{0\leq t\leq T}||u^m||^2_{H^2}+\int^T_0(||\nabla u^m||^2_{H^2}+||\nabla u^m_t||^2_{L^2})\leq C,\\
&|d^m|\leq 1\ \ \rm{a.e.}\Omega\times[0,T],\\
&\sup_{0\leq t\leq T}||\nabla d^m||^2_{H^2}+\int^T_0(||\Delta d^m||^2_{H^2}+||\Delta d^m_t||^2_{L^2})\leq C,\\
&\theta^m\geq \underline{\theta}_0>0\ \ \rm{a.e.}\Omega\times[0,T],\\
&\sup_{0\leq t\leq T}||\theta^m||^2_{H^2}+\int^T_0(||\nabla \theta^m||^2_{H^2}+||\nabla\theta^m_t||^2_{L^2})\leq C,\\
&\sup_{0\leq t\leq T}||p^m||^2_{H^1}+\int^T_0||p||^2_{H^2}\leq C,
\end{align}
which imply that
\begin{align*}
&u^m\rightharpoonup u\quad \textrm{weakly* in }L^\infty(0,T;H^2),\quad u^m\rightharpoonup u\quad \textrm{weakly in }L^2(0,T;H^3),\\
&u^m_t\rightharpoonup u_t\quad \textrm{weakly in }L^2(0,T;H^1);\\
&d^m\rightharpoonup u\quad \textrm{weakly* in }L^\infty(0,T;H^3),\quad d^m\rightharpoonup u\quad \textrm{weakly in }L^2(0,T;H^4),\\
&d^m_t\rightharpoonup d_t\quad \textrm{weakly in }L^2(0,T;H^2);\\
&\theta^m\rightharpoonup \theta\quad \textrm{weakly* in }L^\infty(0,T;H^2),\quad \theta^m\rightharpoonup \theta\quad \textrm{weakly in }L^2(0,T;H^3),\\
&\theta_t^m\rightharpoonup \theta_t\ \textrm{weakly in }L^2(0,T;H^1);\\
&p^m\rightharpoonup p\quad \textrm{weakly* in }L^\infty(0,T;H^1),\quad p^m\rightharpoonup p\quad \textrm{weakly in }L^2(0,T;H^2).
\end{align*}
As to prove the convergence of sequences $(u^m,d^m,\theta^m,p^m)$, we need the following compactness lemma.
\begin{lem}
{\rm(J.Simon 1987\cite{18})} Assume that $X,B$ and $Y$ are three Banach spaces with $X\hookrightarrow\hookrightarrow B\hookrightarrow Y$. Then the following hold true: \\
(i)If $F$ is a bounded subset of $L^p(0,T;X)$£¬$1\leq p<\infty$, and $\frac{\partial F}{\partial t}=\{\frac{\partial f}{\partial t}|f\in F\}$ is bounded in $L^1(0,T;Y)$, then $F$ is relatively compact in $L^p(0,T;B)$ ;\\
(ii)If $F$ is a bounded subset of $L^\infty(0,T;X)$, and $\frac{\partial f}{\partial t}$ is bounded in $L^r(0,T;Y)$ , where $r>1$, then $F$ is relatively compact in $C([0,T];B)$.
\end{lem}

As we know that $H^2 \hookrightarrow\hookrightarrow H^1\hookrightarrow L^2$, then by applying Lemma 5.1, together with (5.1), we can easily conclude that $\{u^m\}^\infty_{m=1}\subset C(0,T;H^1)$ is relatively compact, which means the existence of convergent subsequence, still denoted by $\{u^m\}$, such that $u^m\rightarrow u$ strongly in $C([0,T];H^1)$. Similarly, we have $d^m\rightarrow d$ strongly in $C([0,T];H^2)$ and $\theta^m\rightarrow \theta$ strongly in $C([0,T];H^1)$.

Next we would like to show that, as $m\rightarrow\infty$, each term of system (4.2) will converge to its related term in the original system (1.8) respectively. Here we select only two terms, which are more complex relatively, and the others can be shown in the similar manner.

\begin{lem}
Assume that
\begin{center}
$u^m\rightharpoonup u\ $ weakly in $L^2(0,T;H^3)$, $u^m\rightarrow u\ $ strongly in $C([0,T];H^1)$; \\
$d^m\rightharpoonup d\ $weakly in $L^2(0,T;H^3)$, $d^m\rightarrow d\ $ strongly in $C([0,T];H^2)$; \\
$\theta^m\rightharpoonup \theta\ $weakly in $L^2(0,T;H^3)$, $\theta^m\rightarrow \theta\ $strongly in $C([0,T];H^1)$.
\end{center}
Then
\begin{align}
&\lim_{m\rightarrow\infty}\int^T_0\int\mathrm{div}(\lambda(\theta^m)\nabla d^m\odot\nabla d^m)=\int^T_0\int\mathrm{div}(\lambda(\theta)\nabla d\odot\nabla d),\\
&\lim_{m\rightarrow\infty}\int^T_0\int\mu(\theta^m)(\nabla u^m+\nabla^Tu^m):\nabla u^m=\int^T_0\int\mu(\theta)(\nabla u+\nabla^Tu):\nabla u.
\end{align}
\end{lem}
\proof
We first prove (5.7). Note that
\begin{align}
&\mathrm{div}(\lambda(\theta^m)\nabla d^m\odot\nabla d^m)-\mathrm{div}{\lambda(\theta)\nabla d\odot\nabla d}\nonumber\\
=&\mathrm{div}[(\lambda(\theta^m)-\lambda(\theta))\nabla d^m\odot\nabla d^m]+\mathrm{div}[\lambda(\theta)(\nabla d^m-\nabla d)\odot\nabla d^m]\nonumber\\
&+\mathrm{div}[\lambda(\theta)\nabla d\odot(\nabla d^m-\nabla d)]\nonumber\\
=&\lambda^\prime(\theta^m)\nabla\theta^m\cdot\nabla d^m\odot\nabla d^m-\lambda^\prime(\theta)\nabla\theta\nabla d\odot\nabla d\nonumber\\
&+\lambda(\theta^m)\Delta d^m\cdot\nabla d^m-\lambda(\theta)\Delta d\cdot\nabla d+\lambda(\theta^m)\nabla d^m\cdot\nabla^2 d^m-\lambda(\theta)\nabla d\cdot\nabla^2 d\nonumber\\
=&(\lambda^\prime(\theta^m)-\lambda^\prime(\theta))\nabla \theta^m\cdot\nabla d^m\odot\nabla d^m+\lambda^\prime(\theta)(\nabla\theta^m-\nabla\theta)\nabla d^m\odot\nabla d^m\nonumber\\
&+\lambda^\prime(\theta)\nabla\theta\cdot(\nabla d^m-\nabla d)\odot\nabla d^m+\lambda^\prime(\theta)\nabla\theta\cdot\nabla d\odot(\nabla d^m-\nabla d)\nonumber\\
&+(\lambda(\theta^m)-\lambda(\theta))\Delta d^m\cdot\nabla d^m+\lambda(\theta)(\Delta d^m-\Delta d)\cdot\nabla d^m+\lambda(\theta)\Delta d\cdot(\nabla d^m-\nabla d)\nonumber\\
&+(\lambda(\theta^m)-\lambda(\theta))\nabla d^m\cdot\nabla^2d^m+\lambda(\theta)(\nabla d^m-\nabla d)\cdot\nabla^2d^m+\lambda(\theta)\nabla d\cdot(\nabla^2d^m-\nabla^2d^m).
\end{align}
By applying the differential mean value theorem and directly calculating, we have
\begin{align}
&\mathrm{div}(\lambda(\theta^m)\nabla d^m\odot\nabla d^m)-\mathrm{div}{\lambda(\theta)\nabla d\odot\nabla d}\nonumber\\
\lesssim&|\theta^m-\theta||\nabla d^m|(|\nabla\theta^m||\nabla d^m|+|\nabla^2 d^m|)+|\nabla d^m-\nabla d||\nabla\theta|(|\nabla d^m|+|\nabla d|)\nonumber\\
&+|\nabla d^m-\nabla d|(|\nabla d|+|\nabla d^m|)+(|\nabla^2d^m-\nabla^2d|)(|\nabla d^m|+|\nabla d|)+|\nabla\theta^m-\nabla\theta||\nabla d^m|^2.
\end{align}
Thus,
\begin{align}
&\int^T_0\int|\mathrm{div}[\lambda(\theta^m)\nabla d^m\odot\nabla d^m-\lambda(\theta)\nabla d\odot\nabla d]|\nonumber\\
\lesssim&\int_0^T\int|\theta^m-\theta||\nabla d^m|(|\nabla\theta^m||\nabla d^m|+|\nabla^2 d^m|)\nonumber\\
&+\int^T_0\int|\nabla d^m-\nabla d||\nabla\theta|(|\nabla d^m|+|\nabla d|)+\int_0^T\int|\nabla d^m-\nabla d|(|\nabla d|+|\nabla d^m|)\nonumber\\
&+\int_0^T\int(|\nabla^2d^m-\nabla^2d|)(|\nabla d^m|+|\nabla d|)+\int_0^T\int|\nabla\theta^m-\nabla\theta||\nabla d^m|^2\nonumber\\
\lesssim&\int^T_0||\theta^m-\theta||_{L^2}(||\nabla d^m||^2_{L^\infty}||\nabla \theta^m||_{L^2}+\int^T_0||\nabla d^m||_{L^\infty}||\nabla^2 d^m||_{L^2})\nonumber\\
&+\int_0^T||\nabla d^m-\nabla d||_{L^2}||\nabla \theta||_{L^2}(||\nabla d^m||_{L^\infty}+||\nabla d||_{L^\infty})+\int^T_0||\nabla \theta^m-\nabla\theta||_{L^2}||\nabla d^m||^2_{L^4}\nonumber\\
&+\int_0^T(||\nabla d^m-\nabla d||_{L^2}+||\nabla^2 d^m-\nabla^2d||_{L^2})(||\nabla d^m||_{L^2}+||\nabla d||_{L^2})\nonumber\\
\lesssim &\max_{0\leq t\leq T}||\theta^m-\theta||_{L^2}+\max_{0\leq t\leq T}||\nabla \theta^m-\nabla\theta||_{L^2}\nonumber\\
&+\max_{0\leq t\leq T}||\nabla d^m-\nabla d||_{L^2}+\max_{0\leq t\leq T}||\nabla^2 d^m-\nabla^2d||_{L^2}).
\end{align}
Therefore, the convergence of ${d^m},{\theta^m}$ easily implies
\[\left|\int^T_0\int\mathrm{div}[\lambda(\theta^m)\nabla d^m\odot\nabla d^m-\lambda(\theta)\nabla d\odot\nabla d]\right|\rightarrow 0 (m\rightarrow\infty)\]

Next, we prove (5.8). Similarly, we have
\begin{align}
&\int^T_0\int|\mu(\theta^m)(\nabla u^m+\nabla^Tu^m):\nabla u^m-\mu(\theta)(\nabla u+\nabla^Tu):\nabla u|\nonumber\\
\leq&\int^T_0\int|(\mu(\theta^m)-\mu(\theta))(\nabla u^m+\nabla^Tu^m):\nabla u^m|\nonumber\\
&+\int^T_0\int|\mu(\theta)[\nabla(u^m-u)+\nabla^T(u^m-u)]:\nabla u^m|\nonumber\\
&+\int^T_0\int|\mu(\theta)(\nabla u+\nabla^Tu):(\nabla u^m-\nabla u)|\nonumber\\
\lesssim &\int^T_0\int|\theta^m-\theta||\nabla u^m|^2+\int^T_0\int|\nabla u^m-\nabla u|(|\nabla u^m|+|\nabla u|)\nonumber\\
\lesssim&\int_0^T||\theta^m-\theta||_{L^2}||\nabla u^m||^2_{L^4}+\int^T_0||\nabla u^m-\nabla u||_{L^2}(||\nabla u^m||_{L^2}+||\nabla u||_{L^2})\nonumber\\
\lesssim&\max_{0\leq t\leq T}||\theta^m-\theta||_{L^2}+\max_{0\leq t\leq T}||\nabla u^m-\nabla u||_{L^2}\rightarrow 0\quad(m\rightarrow \infty),
\end{align}
which completes the proof.
\endproof

In conclusion, we have proved that $(u(x,t),\theta(x,t),d(x,t),p(x,t))$ do be strong solutions to system (1.8) over $\Omega\times[0,T]$ and satisfy
\begin{align}
&u(x,t)\in L^\infty(0,T;H^2)\cap L^2(0,T;H^3),u_t\in L^2(0,T;H^1);\nonumber\\
&d(x,t)\in L^\infty(0,T;H^3)\cap L^2(0,T;H^4)\cap L^\infty(\Omega\times[0,T]);\nonumber\\
&d_t(x,t)\in L^2(0,T;H^2),|d(x,t)|\leq 1\ a.e.\Omega\times[0,T];\nonumber\\
&\theta(x,t)\in L^\infty(0,T;H^2)\cap L^2(0,T;H^3);\theta_t(x,t)\in L^2(0,T;H^1);\nonumber\\
&\theta(x,t)\geq\underline{\theta}_0>0\ a.e.\Omega\times[0,T],p(x,t)\in L^\infty(0,T;H^1)\cap L^2(0,T;H^2).\nonumber
\end{align}
\section{Uniqueness of the local strong solutions}
Assume that $(u_i,d_i,\theta_i,p_i),i=1,2$ are two strong solutions to (1.8) with the same initial data (1.9)-(1.10) over $\Omega\times[0,T]$. Then $(u_i,d_i,\theta_i,p_i),(i=1,2)$ satisfy (5.1)-(5.6), and $\bar{u}=u_1-u_2,\bar{d}=d_1-d_2,\bar{\theta}=\theta_1-\theta_2,\bar{p}=p_1-p_2$ solve the following initial boundary problem.
\begin{equation}
\begin{cases}
\bar{u}_t+\bar{u}\cdot\nabla u_1+u_2\cdot\nabla\bar{u}+\nabla\bar{p}=\mathrm{div}[(\mu(\theta_1)-\mu(\theta_2))(\nabla u_1+\nabla^T u_1)]\\
\quad\quad\quad\quad-\mathrm{div}[\mu(\theta_2)(\nabla \bar{u}+\nabla^T \bar{u})]-\mathrm{div}[\lambda(\theta_1)\nabla\bar{d}\odot\nabla d_1]\\
\quad\quad\quad\quad+\mathrm{div}[\lambda(\theta_1)\nabla d_2\odot\nabla\bar{d}+(\lambda(\theta_1)-\lambda(\theta_2))\nabla d_2\odot\nabla d_2]\\
\mathrm{div}\bar{u}=0\\
\bar{\theta}_t+\bar{u}\cdot\nabla\theta_1+u_2\cdot\nabla\bar{\theta}=\Delta \bar{\theta}+(\mu(\theta_1)-\mu(\theta_2))(\nabla u_1+\nabla^T u_1):\nabla u_1\\
\quad\quad\quad\quad+\mu(\theta_2)(\nabla \bar{u}+\nabla^T \bar{u}):\nabla u_1+\mu(\theta_2)(\nabla u_2+\nabla u_2):\nabla\bar{u}\\
\quad\quad\quad\quad-\lambda(\theta_1)\nabla\bar{d}\odot\nabla d_1:\nabla u_1-\lambda(\theta_1)\nabla d_2\odot\nabla\bar{d}:\nabla u_1\\
\quad\quad\quad\quad-\lambda(\theta_1)\nabla d_2\odot\nabla d_2:\nabla\bar{u}-(\lambda(\theta_1)-\lambda(\theta_2))\nabla d_2\odot\nabla d_2:\nabla u_2\\
\bar{d}_t+\bar{u}\cdot\nabla d_1+u_2\cdot\nabla\bar{d}=\Delta\bar{d}-(|d_1|^2-1)\bar{d}-\bar{d}\cdot(d_1+d_2)d_2\\
(\bar{u},\bar{\theta},\bar{d})|_{t=0}=(0,0,0)\\
\bar{u}(x-De_i,t)=\bar{u}(x+De_i,t),\bar{\theta}(x-De_i)=\bar{\theta}(x+De_i),\\
\bar{d}(x-De_i,t)=\bar{d}(x+De_i,t),\forall x\in\partial\Omega,i=1,2,3.
\end{cases}
\end{equation}
What we are going to show is that $(\bar{u},\bar{\theta},\bar{d},\bar{p})=(0,0,0,const.)$ a.e. $[0,T]\times\Omega$.

First, multiplying $(6.1)_1$ by $\bar{u}$ and integrating by parts over $\Omega$, we have
\begin{align}
\frac{1}{2}\frac{\mathrm{d}}{\mathrm{d}t}\int|\bar{u}|^2=&-\int(\mu(\theta_1)-\mu(\theta_2))(\nabla u_1+\nabla^Tu_1):\nabla\bar{u}-\frac{1}{2}\int\mu(\theta_2)|\nabla \bar{u}+\nabla^T \bar{u}|^2\nonumber\\
&+\int\lambda(\theta_1)\nabla\bar{d}\odot\nabla d_1:\nabla\bar{u}+\int\lambda(\theta_1)\nabla d_2\odot\nabla\bar{d}:\nabla\bar{u}\nonumber\\
&+\int(\lambda(\theta_1)-\lambda(\theta_2))\nabla d_2\odot\nabla d_2:\nabla\bar{u}-\int\bar{u}\cdot\nabla u_1\cdot\bar{u},
\end{align}
where the fact $\mathrm{div}u_2=0$ is used.

Second, multiplying $(6.1)_3$ by $\bar{\theta}$ and integrating by parts over $\Omega$, we obtain
\begin{align}
\frac{1}{2}\frac{\mathrm{d}}{\mathrm{d}t}\int|\bar{\theta}|^2&+\int_\Omega|\nabla\bar{\theta}|^2=\int(\mu(\theta_1)-\mu(\theta_2))(\nabla u_1+\nabla^Tu_1):\nabla u_1\bar{\theta}\nonumber\\
&+\int\mu(\theta_2)(\nabla\bar{u}+\nabla^T\bar{u}):\nabla u_1\bar{\theta}+\int\mu(\theta_2)(\nabla u_2+\nabla u_2):\nabla\bar{u}\bar{\theta}\nonumber\\
&-\int\lambda(\theta_1)\nabla\bar{d}\odot\nabla d_1:\nabla u_1\bar{\theta}-\int\lambda(\theta_1)\nabla d_2\odot\nabla\bar{d}:\nabla u_1\bar{\theta}-\int\bar{u}\cdot\nabla\theta_1\bar{\theta}\nonumber\\
&\quad-\int\lambda(\theta_1)\nabla d_2\odot\nabla d_2:\nabla\bar{u}\bar{\theta}-\int[\lambda(\theta_1)-\lambda(\theta_2)]\nabla d_2\odot\nabla d_2:\nabla u_2\bar{\theta}.
\end{align}
Here the fact $\mathrm{div} u_2=0$ have been used.

Third, multiplying $(6.1)_4$ by $\Delta\bar{d}$ and integrating by parts over $\Omega$, we derive
\begin{align}
\frac{1}{2}\frac{\mathrm{d}}{\mathrm{d}t}\int|\nabla\bar{d}|^2+\int|\Delta\bar{d}|^2&=\int[(1-|d_1|^2)\bar{d}+\bar{d}\cdot(d_1+d_2)d_2]\cdot\Delta\bar{d}\nonumber\\
&\quad+\int\bar{u}\cdot\nabla d_1\cdot\Delta\bar{d}+\int u_2\cdot\nabla\bar{d}\cdot\Delta\bar{d}.
\end{align}

Next, summing up $(6.2),(6.3),(6.4)$ implies that
\begin{align}
&\frac{1}{2}\frac{\mathrm{d}}{\mathrm{d}t}\int(|\bar{u}|^2+|\bar{\theta}|^2+|\nabla\bar{d}|^2)+\int(\underline{\mu}|\nabla \bar{u}|^2+|\nabla\bar{\theta}|^2+|\Delta\bar{d}|^2)\nonumber\\
\leq&-\int(\mu(\theta_1)-\mu(\theta_2))(\nabla u_1+\nabla^Tu_1):\nabla\bar{u}+\int\lambda(\theta_1)\nabla\bar{d}\odot\nabla d_1:\nabla\bar{u}\nonumber\\
&+\int\lambda(\theta_1)\nabla d_2\odot\nabla\bar{d}:\nabla\bar{u}+\int(\lambda(\theta_1)-\lambda(\theta_2))\nabla d_2\odot\nabla d_2:\nabla\bar{u}\nonumber\\
&-\int\bar{u}\cdot\nabla u_1\cdot\bar{u}+\int(\mu(\theta_1)-\mu(\theta_2))(\nabla u_1+\nabla^Tu_1):\nabla u_1\bar{\theta}\nonumber\\
&+\int\mu(\theta_2)(\nabla\bar{u}+\nabla^T\bar{u}):\nabla u_1\bar{\theta}+\int\mu(\theta_2)(\nabla u_2+\nabla u_2):\nabla\bar{u}\bar{\theta}\nonumber\\
&-\int\lambda(\theta_1)\nabla\bar{d}\odot\nabla d_1:\nabla u_1\bar{\theta}-\int\lambda(\theta_1)\nabla d_2\odot\nabla\bar{d}:\nabla u_1\bar{\theta}-\int\bar{u}\cdot\nabla\theta_1\bar{\theta}\nonumber\\
&-\int\lambda(\theta_1)\nabla d_2\odot\nabla d_2:\nabla\bar{u}\bar{\theta}-\int[\lambda(\theta_1)-\lambda(\theta_2)]\nabla d_2\odot\nabla d_2:\nabla u_2\bar{\theta}\nonumber\\
&+\int[(1-|d_1|^2)\bar{d}+\bar{d}\cdot(d_1+d_2)d_2]\cdot\Delta\bar{d}+\int\bar{u}\cdot\nabla d_1\cdot\Delta\bar{d}+\int u_2\cdot\nabla\bar{d}\cdot\Delta\bar{d}\nonumber\\
=&\sum^{16}_{i=1}K_i,
\end{align}
where we have used $\int\underline{\mu}|\nabla\bar{u}|^2\leq\int\frac{\mu(\theta_1)}{2}|\nabla\bar{u}+\nabla^T\bar{u}|^2$.

As to the desired uniqueness, we will estimate the right-hand side of (6.5) term by term by applying (5.1)-(5.6). We should remind readers that we will always just dealing with the highest order terms while applying interpolation inequality. We do as follows.
\begin{align}
K_1&\lesssim\int|\bar{\theta}||\nabla u_1||\nabla\bar{u}|\lesssim||\bar{\theta}||_{L^2}||\nabla u_1||_{L^\infty}||\nabla\bar{u}||_{L^2}\nonumber\\
&\lesssim||\bar{\theta}||_{L^2}||\nabla\Delta u_1||_{L^2}||\nabla\bar{u}||_{L^2}\nonumber\\
&\leq\frac{\underline{\mu}}{8}||\nabla\bar{u}||^2_{L^2}+C||\nabla\Delta u_1||^2_{L^2}||\bar{\theta}||^2_{L^2}.\\
K_2+K_3&\lesssim\int(|\nabla\bar{d}||\nabla\bar{u}|(|\nabla d_1|+|\nabla d_2|)\nonumber\\
&\lesssim||\nabla\bar{d}||_{L^2}||\nabla\bar{u}||_{L^2}(||\nabla d_1||_{L^\infty}+||\nabla d_2||_{L^\infty})\nonumber\\
&\lesssim||\nabla\bar{d}||_{L^2}||\nabla\bar{u}||_{L^2}(||\nabla\Delta d_1||_{L^2}+||\nabla\Delta d_2||_{L^2})\nonumber\\
&\lesssim||\nabla\bar{d}||_{L^2}||\nabla\bar{u}||_{L^2}\nonumber\\
&\leq\frac{\underline{\mu}}{8}||\nabla\bar{u}||^2_{L^2}+C||\nabla\bar{d}||^2_{L^2}.\\
K_4&\lesssim\int|\bar{\theta}||\nabla d_2|^2|\nabla \bar{u}|\lesssim||\nabla d_2||^2_{L^\infty}||\bar{\theta}||_{L^2}||\nabla\bar{u}||_{L^2}\nonumber\\
&\lesssim||\nabla\Delta d_2||^2_{L^2}||\bar{\theta}||_{L^2}||\nabla\bar{u}||_{L^2}\lesssim||\bar{\theta}||_{L^2}||\nabla\bar{u}||_{L^2}\nonumber\\
&\leq\frac{\underline{\mu}}{8}||\nabla\bar{u}||^2_{L^2}+C||\bar{\theta}||^2_{L^2}.\\
K_5&\lesssim\int|\bar{u}|^2|\nabla u_1|\lesssim||\bar{u}||^2_{L^4}||\nabla u_1||_{L^2}\lesssim||\bar{u}||^{\frac{1}{2}}_{L^2}||\nabla\bar{u}||^{\frac{3}{4}}_{L^2}\nonumber\\
&\leq\frac{\underline{\mu}}{8}||\nabla\bar{u}||^2_{L^2}+C||\bar{u}||^2_{L^2}.\\
K_6&\lesssim\int|\bar{\theta}|^2|\nabla u_1|^2\lesssim||\bar{\theta}||^2_{L^2}||\nabla u_1||^2_{L^\infty}\leq C||\nabla\Delta u_1||^2_{L^2}||\bar{\theta}||^2_{L^2}\\
K_7+K_8&\lesssim\int(|\nabla u_1|+|\nabla u_2|)|\nabla\bar{u}||\bar{\theta}|\nonumber\\
&\lesssim(||\nabla u_1||_{L^\infty}+||\nabla u_2||_{L^\infty})||\nabla\bar{u}||_{L^2}||\bar{\theta}||_{L^2}\nonumber\\
&\lesssim(||\nabla\Delta u_1||_{L^2}+||\nabla\Delta u_2||_{L^2})||\nabla\bar{u}||_{L^2}||\bar{\theta}||_{L^2}\nonumber\\
&\leq \frac{\underline{\mu}}{8}||\nabla\bar{u}||^2_{L^2}+C(||\nabla\Delta u_1||^2_{L^2}+||\nabla\Delta u_2||^2_{L^2})||\bar{\theta}||^2_{L^2}\\
K_9+K_{10}&\lesssim\int(|\nabla d_1|+|\nabla d_2|)|\nabla u_1||\nabla\bar{d}||\bar{\theta}|\nonumber\\
&\lesssim(||\nabla d_1||_{L^\infty}+||\nabla d_2||_{L^\infty})||\nabla u_1||_{L^3}||\nabla\bar{d}||_{L^6}||\bar{\theta}||_{L^2}\nonumber\\
&\lesssim(||\nabla\Delta d_1||_{L^2}+||\nabla\Delta d_2||_{L^2})||\Delta u_1||_{L^2}||\Delta\bar{d}||_{L^2}||\bar{\theta}||_{L^2}\nonumber\\
&\leq\frac{1}{6}||\Delta\bar{d}||^2_{L^2}+C||\bar{\theta}||^2_{L^2}.\\
K_{11}&\lesssim\int|\bar{u}||\nabla\theta_1||\bar{\theta}|\lesssim||\bar{u}||_{L^6}||\nabla\theta_1||_{L^3}||\bar{\theta}||_{L^2}
\lesssim||\Delta\theta_1||_{L^2}||\nabla\bar{u}||_{L^2}||\bar{\theta}||_{L^2}\nonumber\\
&\leq\frac{\underline{\mu}}{8}||\nabla\bar{u}||^2_{L^2}+C||\bar{\theta}||^2_{L^2}.\\
K_{12}+K_{13}&\lesssim\int|\nabla d_2|^2|\nabla\bar{u}||\bar{\theta}|+\int|\nabla d_2|^2|\bar{\theta}|^2|\nabla u_2|\nonumber\\
&\lesssim||\nabla d_2||^2_{L^\infty}||\bar{\theta}||_{L^2}||\nabla\bar{u}||_{L^2}+||\nabla d_2||^2_{L^\infty}||\bar{\theta}||^2_{L^2}||\nabla u_2||_{L^\infty}\nonumber\\
&\lesssim||\nabla\Delta d_2||^2_{L^2}||\bar{\theta}||_{L^2}||\nabla\bar{u}||_{L^2}+||\nabla\Delta d_2||^2_{L^2}||\bar{\theta}||^2_{L^2}||\nabla\Delta u_2||_{L^2}\nonumber\\
&\leq\frac{\underline{\mu}}{8}||\nabla\bar{u}||^2_{L^2}+C||\bar{\theta}||^2_{L^2}+C||\nabla\Delta u_2||^2_{L^2}||\bar{\theta}||^2_{L^2}.\\
K_{14}&=\int[2(\nabla d_1\cdot d_1)\otimes\bar{d}+(|d_1|^2-1)\nabla\bar{d}+\nabla\bar{d}\cdot(d_1+d_2)\otimes d_2]:\nabla\bar{d}\nonumber\\
&\quad+\int[\bar{d}\cdot(\nabla d_1+\nabla d_2)\otimes d_2+\bar{d}\cdot(d_1+d_2)\nabla d_2]:\nabla\bar{d}\nonumber\\
&\lesssim\int(|\nabla d_1||\bar{d}|+|\nabla\bar{d}|+|\bar{d}||\nabla d_2|)|\nabla\bar{d}|\nonumber\\
&\lesssim||\nabla d_1||_{L^2}||\bar{d}||_{L^\infty}||\nabla\bar{d}||_{L^2}+||\nabla\bar{d}||^2_{L^2}+||\bar{d}||_{L^\infty}||\nabla d_2||_{L^2}||\nabla\bar{d}||_{L^2}\nonumber\\
&\lesssim||\Delta\bar{d}||_{L^2}||\nabla\bar{d}||_{L^2}+||\nabla\bar{d}||^2_{L^2}\nonumber\\
&\leq\frac{1}{6}||\Delta\bar{d}||^2_{L^2}+C||\nabla\bar{d}||^2_{L^2}.\\
K_{15}+K_{16}&\lesssim\int|\bar{u}||\nabla d_1||\Delta\bar{d}|+\int|u_2||\nabla\bar{d}||\Delta\bar{d}|\nonumber\\
&\lesssim||\nabla d_1||^2_{L^\infty}||\bar{u}||_{L^2}||\Delta\bar{d}||_{L^2}+||u_2||_{L^\infty}||\nabla\bar{d}||_{L^2}||\Delta\bar{d}||_{L^2}\nonumber\\
&\lesssim||\nabla\Delta d_1||_{L^2}||\bar{u}||_{L^2}||\Delta\bar{d}||_{L^2}+||\Delta u_2||_{L^2}||\nabla\bar{d}||_{L^2}||\Delta\bar{d}||_{L^2}\nonumber\\
&\leq\frac{1}{6}||\Delta\bar{d}||^2_{L^2}+C||\bar{u}||^2_{L^2}+C||\nabla\bar{d}||^2_{L^2}.
\end{align}

Substituting $K_1,\cdots, K_{16}$ into $(6.5)$, we have
\begin{align}
&\frac{1}{2}\frac{\mathrm{d}}{\mathrm{d}t}(||\bar{u}||^2_{L^2}+||\nabla\bar{d}||^2_{L^2}+||\bar{\theta}||^2_{L^2})+(\underline{\mu}||\nabla\bar{u}||^2_{L^2}+||\Delta\bar{d}||^2_{L^2}
+||\nabla\bar{\theta}||^2_{L^2})\nonumber\\
\leq&C(||\bar{u}||^2_{L^2}+||\nabla\bar{d}||^2_{L^2}+||\bar{\theta}||^2_{L^2})+C(||\nabla\Delta u_1||^2_{L^2}+||\nabla\Delta u_2||^2_{L^2})||\bar{\theta}||^2_{L^2}\nonumber\\
&+\frac{1}{2}||\Delta\bar{d}||^2_{L^2}+\frac{7\mu}{8}||\nabla\bar{u}||^2_{L^2},
\end{align}
which is equivalent to
\begin{align}
&\frac{\mathrm{d}}{\mathrm{d}t}(||\bar{u}||^2_{L^2}+||\nabla\bar{d}||^2_{L^2}+||\bar{\theta}||^2_{L^2})+(\underline{\mu}||\nabla\bar{u}||^2_{L^2}+||\Delta\bar{d}||^2_{L^2}+||\nabla\bar{\theta}||^2_{L^2})\nonumber\\
\leq&C(||\nabla\Delta u_1||^2_{L^2}+||\nabla\Delta u_2||^2_{L^2}+1)(||\bar{u}||^2_{L^2}+||\nabla\bar{d}||^2_{L^2}+||\bar{\theta}||^2_{L^2}).
\end{align}
Since we have
\begin{align}
\int^T_0(||\nabla\Delta u_1||^2_{L^2}+||\nabla\Delta u_2||^2_{L^2}+1)\leq C,
\end{align}
the Gronwall inequality implies that
\begin{align}
\sup_{0\leq t\leq T}(||\bar{u}||^2_{L^2}+||\nabla\bar{d}||^2_{L^2}+||\bar{\theta}||^2_{L^2})+\int^T_0(\underline{\mu}||\nabla\bar{u}||^2_{L^2}+||\Delta\bar{d}||^2_{L^2}+||\nabla\bar{\theta}||^2_{L^2})\leq 0,
\end{align}
where we have used $(\bar{u},\bar{\theta},\bar{d})|_{t=0}=(0,0,0)$.

This directly implies $\bar{u}=\nabla\bar{d}=\bar{\theta}=0,a.e\ \Omega\times[0,T]$, and substituting which into $(6.1)_4$, one obtains the equation
\begin{align}
\bar{d}_t=-(|d_1|^2-1)\bar{d}-\bar{d}\cdot(d_1+d_2)d_2.
\end{align}
Multiplying this simplified equation by $\bar{d}$ and integrating over $\Omega$, together with the fact that $0\leq|d_1|\leq 1,0\leq|d_2|\leq 1$, we have
\begin{align}
\frac{1}{2}\frac{\mathrm{d}}{\mathrm{d}t}||\bar{d}||^2_{L^2}\leq C||\bar{d}||^2_{L^2}.
\end{align}
Once again we apply Gronwall inequality and obtain $||\bar{d}||^2_{L^2}\leq 0$, which means $\bar{d}=0,a.e\ \Omega\times[0,T]$.
Finally, substituting $(\bar{u},\bar{d},\bar{\theta})=(0,0,0),a.e\ \Omega\times[0,T]$ into $(6.1)_1$, we have $\nabla\bar{p}=0,a.e\ \Omega\times[0,T]$.

So far, we finish the proof of uniqueness of the strong solution to initial boundary problem $(1.8)-(1.10)$. Of course, $p$ is uniquely determined up to a constant.

\bigbreak\noindent{ \large\bf Acknowledgment} This work is supported by the National
Basic Research Program of China (973 Program) (No.2011CB808002), the National Natural
Science Foundation of China (No.11071086, No.11371152 and No.11128102), the Natural Science Foundation of Guangdong Province (No.S2012010010408), the University Special Research Foundation for Ph.D Program (No.20104407110002) and the Scientific Research Foundation of Graduate School of South China Normal University (No.2014ssxm04).



\end{document}